\documentclass[preprint,12pt,3p]{elsarticle}
\usepackage{lineno,hyperref}
\usepackage{graphics,subfigure,varwidth}
\usepackage{caption,booktabs}
\usepackage{bm}
\usepackage{amsmath}
\usepackage{float,epstopdf}
\usepackage{color}
\usepackage{multirow}

\def\rma{\mathrm{Ma}}

\def\rmd{\mathrm{d}}
\def\rme{\mathrm{e}}
\def\rmQ{\bm{\mathrm Q}}
\def\rmU{\bm{\mathrm U}}
\def\rmF{\bm{\mathrm F}}
\def\rmx{\bm{\mathrm x}}
\def\rmu{\bm{\mathrm u}}

\def\rmxi{\bm{\mathrm \xi}}
\def\rmpsi{\bm{\mathrm \psi}}

\modulolinenumbers[5]

\begin{document}

\begin{frontmatter}

\title{A two-stage fourth-order gas-kinetic CPR method for the Navier-Stokes equations on triangular meshes}

\author[address1]{Chao Zhang}
\ead{zhangcha15@mails.tsinghua.edu.cn}
\author[address1]{Qibing Li\corref{mycorrespondingauthor}}
\cortext[mycorrespondingauthor]{Corresponding author}
\ead{lqb@tsinghua.edu.cn}
\author[address2]{Z.J.Wang}
\ead{zjw@ku.edu}
\author[address3]{Jiequan Li}
\ead{li_jiequan@iapcm.ac.cn}
\author[address1]{Song Fu}
\ead{fs-dem@tsinghua.edu.cn}
\address[address1]{AML, Department of Engineering Mechanics, Tsinghua University, Beijing 100084, China}
\address[address2]{Department of Aerospace Engineering, The University of Kansas, Lawrence, KS 66045, USA}
\address[address3]{Laboratory of Computational Physics, Institute of Applied Physics and Computational Mathematics, Beijing 100088, China}

\begin{abstract}
A highly efficient gas-kinetic scheme with fourth-order accuracy in both space and time  is developed for the Navier-Stokes equations on triangular meshes. The scheme combines an efficient correction procedure via reconstruction (CPR) framework with a robust gas-kinetic flux formula, which computes both the flux and its time-derivative. The availability of the flux time-derivative makes it straightforward to adopt an efficient two-stage temporal discretization to achieve fourth-order time accuracy. 
In addition, through the gas-kinetic evolution model, the inviscid and viscous fluxes are coupled and computed uniformly without any separate treatment for the viscous fluxes. As a result, the current scheme is more efficient than traditional explicit CPR methods with a separate treatment for viscous fluxes, and a fourth order Runge-Kutta approach. Furthermore, a robust and accurate subcell finite volume (SCFV) limiting procedure is extended to the CPR framework for troubled cells, resulting in subcell resolution of flow discontinuities. Numerical tests demonstrate the high accuracy, efficiency and robustness of the current scheme in a wide range of inviscid and viscous flow problems from subsonic to supersonic speeds.

\end{abstract}


\end{frontmatter}

\section{Introduction}
First and second order methods have been widely used in industrial computational fluid dynamics (CFD) applications due to their simplicity and robustness. With the same computational cost, high-order methods are capable of generating solutions with higher accuracy and thus efficiency\cite{ZJWang2013}. In particular, high-order methods on unstructured meshes are more flexible for complex geometries than their counterparts on structured meshes, and the corresponding workload for mesh generation can be significantly reduced. 
As a result, various high-order methods have been developed in recent decades, which roughly fall into two categories. The first one is the finite volume (FV) method which has good robustness and resolution of discontinuities in high-speed flows \cite{TJBarth1990,CHu1999}. However, the compactness is always a bottleneck. To overcome this problem, many compact high-order FV methods have been developed, \textcolor{black}{in which only face-neighboring cells are involved in reconstruction}, such as the variational reconstruction method \cite{QWang2017}, the explicit multi-step reconstruction method \cite{YZhang2019} and the subcell finite volume method (SCFV) \cite{JPan2017}. In the SCFV method, each element is partitioned into subcells. The reconstruction on each subcell can involve subcells of face-neighbor elements. The residual at each subcell only depends on the current element and its nearest face neighbor elements, thus achieving good compactness while preserving the advantages of a finite volume framework such as the strong robustness and good resolution of discontinuities. The accuracy and efficiency can also be improved in smooth flow regions due to the packed continuous reconstruction for sub cells, \textcolor{black}{in which the solution polynomial is continuous across sub cells.} 

The second one includes discontiuous  methods with internal degrees-of-freedom (IDOF) inside each element. They are compact and can preserve high accuracy in smooth flow regions, such as the Discontinuous Galerkin (DG) method \cite{BCockburn1998,HLuo2008,XZhang2010}. 
In 2007, Huynh developed a high-order method named flux reconstruction (FR) with the advantage of simplicity and efficiency. It solves the 1D conservation laws in differential form \cite{Huynh2007, Huynh2009}. Based on a set of solution points (SPs) as IDOFs, the method evaluates the derivative of the flux by constructing a discontinuous piece-wise flux polynomial. Then a flux correction is added to the discontinuous flux by removing the flux jumps at cell interfaces. The approach provides a unified framework for many existing high-order methods with appropriate choices of correction function, such as the DG, the spectral difference (SD) method \cite{YLiu2006,CLiang2009} and the spectral volume (SV) method \cite{ZJWang2002,ZJWang2004}.  And it has been shown to be simpler and more efficient than the original versions of these methods. Compared to the modal DG method, the cost of the numerical integration can be avoided, and by choosing the SPs properly, the reconstruction can be very efficient. The FR method can be extended to quadrilateral mesh via tensor product directly. And it was first extended to triangular and mixed meshes by Wang and Gao, named lifting collocation penalty (LCP) \cite{ZJWang2009,HGao2009}. Extension to the 3D Euler
and Navier-Stokes (N-S) equations on mixed meshes can be referred to Refs.\cite{THaga2010,THaga2011}. Due to the tight connection between FR and LCP, they are renamed as correction procedure via reconstruction (CPR). Detailed reviews of the CPR method can be found in Ref.\cite{HTHuynh2014}.

A challenging problem for high-order methods based on IDOFs is shock capturing, especially for strong shock waves, as the solution variables are always represented by continuous polynomials within each element. To enhance the resolution of flow discontinuities, a subcell finite volume (SCFV) limiting procedure has been proposed in the DG method and has shown good robustness and accuracy \cite{MDumbser2014,MDumbser2016}. This limiter is applied only to "troubled cells" with the regular DG method applied elsewhere. Thus the high-order accuracy of DG is preserved and the resolution of discontinuities in high-speed flow is enhanced. 

The aforementioned methods mainly focus on the high-order discretization of space. However, high-order time evolution is also very important. Traditional high-order methods usually adopt Riemann solvers to compute the inviscid flux. For viscous flows, the viscous flux needs to be treated additionally. To achieve high-order time accuracy, multi-stage Runge-Kutta (R-K) methods are adopted for time integration. 

Based on the mesoscopic Bhatnagar-Gross-Krook (BGK) model, the gas-kinetic scheme (GKS) offers an alternative way to recover the N-S solutions \cite{KXu2001,QBLi2005}, which is achieved through the first-order Chapman-Enskog expansion. With the help of the local integral solution of the BGK equation, a time-dependent flux function can be constructed by a Taylor expansion of the gas distribution function in space and time. Different from traditional Riemann solvers, the inviscid and viscous fluxes in GKS are coupled and obtained simultaneously. It is genuinely multi-dimensional by involving both normal and tangential derivatives in the gas distribution function \cite{QBLi2006}. The multiscale evolution process from a kinetic scale to a hydrodynamic scale keeps a good balance between accuracy and robustness \cite{KXu2005}.

Through a second-order Taylor expansion of the gas distribution function, a third-order multi-dimensional GKS (HGKS or HBGK) approach has been developed successfully \cite{QBLi2008,QBLi2010,QBLi2012}. Based on this high-order gas evolution model, a series of third-order gas-kinetic schemes have been developed within a single stage, such as the compact GKS \cite{LPan2015,LPan2016}, SV-GKS \cite{NLiu2017}, DG-GKS \cite{XDRen2015} CLS-GKS \cite{JL2019}, SCFV-GKS \cite{CZhang2019}, etc. Recently, a single-stage third-order gas-kinetic CPR method on triangular meshes has also been developed \cite{CZhang2018}. By combining the efficient CPR framework with the third-order gas-kinetic flux, it shows high accuracy and efficiency in many benchmark flow problems. Therefore, it is worth developing a fourth-order gas-kinetic CPR scheme on triangular meshes to achieve higher accuracy and efficiency. 

A straightforward way to develop a fourth-order gas-kinetic CPR is to use a third-order Taylor expansion in space and time to construct a single stage time-evolution flux. Through this way a one-dimensional finite volume scheme has been developed \cite{NLiu2014}. Fortunately, the two-stage fourth-order time-stepping method has been developed for Lax-Wendroff type flow solvers \cite{JLi2016, ZDu2018, JQLi201902} by using the flux and its first-order time derivative   which seems simpler. Recently, the two-stage temporal discretization has also been extended to DG based on the flux solver for the generalized Riemann problem (GRP) for the Euler equations  \cite{JQLi2019}.  As only one middle stage is used and little additional computational cost for the time derivative is required, the scheme shows higher efficiency than the traditional RKDG. In addition, the reduction of temporal stage also means the reduction of effective stencil size, or better locality, which benefits parallel computations. In fact, the inherent space-time coupling in this two-stage technique, as well as the locality is important for Euler and compressible N-S equations. Since the gas-kinetic flux also provides a time-evolving flux function, an efficient two-stage fourth-order GKS has also been constructed successfully through the combination of the second-order gas-kinetic flux solver and the two-stage temporal discretization method \cite{LPan2016II, FXZhao2019}. 
In the meanwhile, the multi-stage multi-derivative method was developed for hyperbolic conservation laws \cite{DCSeal2014}, which has also been used to develop a family of HGKS \cite{XJi2018}.

The objective of current study is to construct a highly efficient forth-order scheme on triangular meshes for the compressible N-S equations. To achieve high efficiency, the new scheme fully combines the efficient CPR framework with the efficient two-stage fourth-order time stepping method, as well as the robust time-evolving second-order gas-kinetic flux solver. Furthermore, a robust SCFV limiting approach is extended to CPR to enhance the resolution of discontinuities. Compared with existing two-stage fourth-order finite volume GKS, the present method is suitable for triangular meshes, and it is compact and efficient. On the other hand, compared with existing CPR methods, it can achieve higher efficiency due to the coupling of inviscid and viscous effects in the gas-kinetic flux, and the efficient time stepping technique with only one middle stage instead of the multi-stage R-K method. In addition, 
with the help of SCFV, the limiting procedure is simple and the subcell resolution of discontinuities can be automatically achieved in high-speed flows.

The paper is organized as follows. In Section \ref{section2}, the construction of the current method is presented, including the CPR framework, the gas-kinetic flux solver, the two-stage 4th order temporal discretization and the subcell finite volume limiting procedure. Numerical tests are presented in Section \ref{section3} to verify and demonstrate the performance of the current scheme. The last section draws the conclusions.

\section{Numerical method} \label{section2}
\subsection{CPR framework}
\begin{figure}[H]
  \centering
  \includegraphics[scale=0.3]{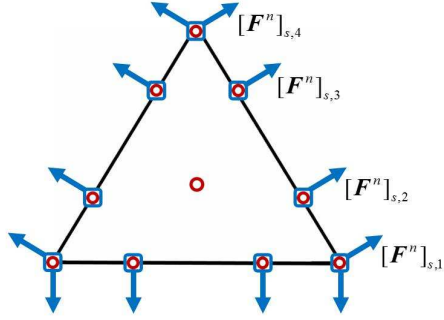}\\
  \captionsetup{labelsep=period}
  \caption{Solution points (circles) and flux points (squares) with $k$=3.}
  \label{SPsFPs}
\end{figure}

The CPR framework has many competitive features compared to other high-order methods. 
Compact reconstruction under the CPR framework is straightforward. It provides a unified framework for DG, SD, SV, etc, and is more efficient and simpler to be implemented. 
Here the CPR framework is briefly reviewed. Consider the conservation law,
\begin{equation}\label{eq_ConLaw}
    \frac{\partial\rmQ}{\partial{t}}+\nabla\cdot\rmF =0,
\end{equation}
where $\rmQ=(\rho,\rho\rmU,\rho E)^T$ are the conservative variables, in which $\rmU=(U,V)$ are the macroscopic velocities and $\rmF=(\rmF_x,\rmF_y)$ is the flux vector. The computational domain is divided into $N$ non-overlapping triangular cells $\{\Omega_i\}$. In each cell ${\Omega}_i$, $\rmQ$ is approximated by a solution polynomial $\rmQ_i$ with degree $k$, which belongs to the polynomial space $P^k$. $\rmQ_i$ can be constructed by the Lagrange interpolation,
\begin{equation}\label{eq_Q(x,y)}
   \rmQ_i(\rmx)=\sum_j L_j(\rmx)\rmQ_{i,j},
\end{equation}
where ${L}_j(\rmx)$ is the Lagrange basis, $\rmQ_{i,j}$ are the conservative variables at the solution point (SP) with the position $\rmx_{i,j}=(x_{i,j},y_{i,j})$. For the fourth-order CPR, $k=3$, the number of SPs is equal to $(k+1)(k+2)/2=10$. To involve the interaction between cells, we need to define $(k+1)$ flux points (FPs) at each cell interface to compute the common flux. The distribution of SPs and FPs is shown in Fig.\ref{SPsFPs}. For efficiency, the SPs are chosen to coincide with the FPs.
The semi-discrete CPR formulation can be expressed as
\begin{equation}\label{eq_CPR}
\begin{split}
    \frac{\partial\rmQ_{i,j}}{\partial{t}}&= \mathcal{R}_{i,j}({\rmF}), \\  
\mathcal{R}_{i,j}({\rmF)}=-\Pi_j\left(\nabla\cdot{\rmF_i}\right)&-\frac{1}{|{\Omega}_i|}\sum_{s\in\partial{\Omega}_i}\sum_{l}\alpha_{j,s,l}[{\rmF}_{i,i+}]_{s,l}|\Gamma|_s,
\end{split}
\end{equation}
where the second term denotes the flux divergence at SPs, $\Pi$ is the projection operator which projects the flux divergence term onto $P^k$. Through the Lagrange interpolation, a flux polynomial ${\rmF}_i={\rmF}(\rmQ_i)$ with degree $k$ can also be constructed. Then the projection can be expressed as
\begin{equation}\label{eq4}
   \Pi{\left(\nabla\cdot{\rmF}_i\right)}=\nabla\cdot\left(\sum_{j}{L}_j(\rmx){\rmF}(\rmQ_{i,j})\right),
\end{equation}
where ${\rmF}(\rmQ_{i,j})$ is the flux at SPs. The third term in Eq.(\ref{eq_CPR}), named correction field, is to involve the interaction between cell ${\Omega}_i$ and its face neighbors ${\Omega}_{i+}$, in which $[{\rmF}_{i,i+}]_{s,l}=[{\rmF}_{com}(\rmQ_i,\rmQ_{i+},\bm{n})-{\rmF}(\rmQ_i)\cdot\bm{n}]_{s,l}$  is the normal flux difference at FPs and ${\rmF}_{com}(\rmQ_i,\rmQ_{i+},\bm{n})$ is the common flux. $\alpha_{j,s,l}$  is the lifting coefficient, which is independent of the solution and geometry. $|{\Omega}_i|$ is the area of cell $i$  and $|\Gamma|_s$  is the length of triangular edge $s$. More details of the CPR framework can be found in Ref.\cite{Huynh2007, ZJWang2009}.

\subsection{Gas-kinetic flux solver}\label{section2.3}
Traditional Riemann solvers are usually used to compute the fluxes at the flux points. Then the governing equation is discretized in time using a multi-stage R-K scheme. In order to achieve 4th order accuracy in time, a minimum number of 4 stages is needed in the R-K scheme. In the present study, we improve the efficiency of the CPR method combining the gas-kinetic flux solver and a two-stage fourth-order time marching scheme. In particular, the second-order gas-kinetic flux solver is the foundation of the overall algorithm. Its basic principle is briefly introduced here.

In the mesoscopic gas-kinetic theory, the flow is described by the gas distribution function $f=f(\rmx,t,\rmu,\rmxi)$ which is a function of physical space $\rmx$, time $t$, particle velocity $\rmu=(u,v)$ and internal degrees of freedom $\rmxi$. In the present study, only two-dimensional (2D) flow is considered.
The macroscopic conservative variables $\rmQ$ and the flux vector $\rmF$ can be obtained by taking moments of $f$ in the phase space, i.e.,
\begin{equation}\label{eq_f_Q}
  \rmQ=\int f\rmpsi\rmd\Xi,~~~
\end{equation}
\begin{equation}\label{eq_f_F}
  {\rmF}_\sigma=\int u_\sigma f\rmpsi\rmd\Xi,~ \sigma=1,2,
\end{equation}
where $\rmpsi={\left(1,\rmu,({\rmu}^2+{\rmxi}^2)/2\right)}^T$ is the vector of moments, $\rmd\Xi=\rmd\rmu\rmd\rmxi$ is the element of the phase space. The governing equation of $f$ is the 2D BGK equation \cite{Bhatnagar1954}
\begin{equation}\label{eq_bgk}
\frac{\partial f}{\partial t}+\rmu\cdot\nabla{f}= \frac{g-f}{\tau},
\end{equation}
where $\tau=\mu/p$ is the collision time dependent on the dynamic viscous coefficient  $\mu$ and pressure $p$. The local equilibrium state $g$ approached by $f$ is the Maxwellian distribution,
 \begin{equation}
 g=\rho (2\pi RT)^{-(K+2)/2} \rme^{-[ (\rmu-\rmU)^2+\rmxi^2]/(2RT)},
\end{equation}
where $\rho$ is the density, $R$ is the gas constant, $T$ is the temperature, and $K$ is the total number of $\rmxi$, which is equal to $(5-3\gamma)/(\gamma-1)+1$ for two-dimensional flow, in which $\gamma$ is the specific heat ratio. By taking moments of Eq.(\ref{eq_bgk}), the macroscopic conservation law Eq.(\ref{eq_ConLaw}) can be recovered, in which the collision term $(g-f)/\tau$ vanishes automatically due to the conservation of mass, moments and total energy during collisions, i.e., the compatibility condition. Particularly, the Naiver-Stokes equations can be recovered through the first-order Chapman-Enskog expansion \cite{KXu2015}
\begin{equation}\label{eq_C_E}
\begin{split}
f_{NS}=g-\tau \left(\frac{\partial g}{\partial t}+\rmu\cdot\nabla{g}\right).
\end{split}
\end{equation}

To update the numerical solution under the CPR framework, the fluxes at FP and SP need to be determined. Since the solution is discontinuous across cell interfaces, common fluxes are computed with the solutions from both sides of the cell interface. In the present study, the time-dependent gas distribution function is constructed based on the local analytical solution of the BGK equation, i.e.,
\begin{equation}\label{eq_BGK_exact}
\begin{split}
 f(\rmx,t,\rmu,\rmxi)=\frac{1}{\tau}\int_0^t g(\rmx -\rmu (t-t'),t',\rmu ,\rmxi )\rme^{-(t-t')/\tau}\rmd t'
                +\rme^{-t/\tau}f_0(\rmx -\rmu t,\rmu ,\rmxi),
\end{split}
\end{equation}
where $f_0$ is the piecewise continuous initial distribution function at the start of each time step, $g$ is the local equilibrium state. For simplicity, the cell interface is assumed perpendicular to the $x$-axis. Through a first-order Taylor expansion of $f_{NS}$ and $g$ at a FP to construct $f_0$ and $g$ in the neighborhood, the time-dependent gas distribution function can be obtained with Eq.(\ref{eq_BGK_exact}), i.e.,
\begin{equation}\label{eq_f_dis}
\begin{split}
 f_{_\mathrm{FP}}(t,\rmu,\rmxi)=&g_0\left(1-\rme^{-t/\tau}
                 +((t+\tau)\rme^{-t/\tau}-\tau)(a_1u+a_2v)
                 +(t-\tau+\tau\rme^{-t/\tau})A\right) \\
        &+\mathrm{e}^{-t/\tau}g_R\left(1-(\tau+t)(a_1^Ru+a_2^Rv)-\tau A^R\right) \mathrm{H}(u) \\
        &+\mathrm{e}^{-t/\tau}g_L\left(1-(\tau+t)(a_1^Lu+a_2^Lv)-\tau A^L\right) \left(1-\mathrm{H}(u)\right),
\end{split}
\end{equation}
where $\mathrm{H}(u)$ is the Heaviside function. The coefficients $a_1,a_2$ and $A$ are related to the Taylor expansion of the corresponding Maxwellian functions, i.e., the derivatives of $g_0$. Other coefficients with superscript $L, R$ correspond to $g_L$ and $g_R$ respectively. More details can be found in Refs.\cite{KXu2001,QBLi2005}. Both normal and tangential spatial derivatives are involved in the construction of the gas distribution function, which depicts a multidimensional transport process across the interface, contributing to a intrinsically multidimensional flux solver\cite{QBLi2006}. 

Since the solution is continuous inside each cell, the gas distribution function at a SP can be simplified as,
\begin{equation}\label{eq_f_con}
\begin{split}
  f_{_\mathrm{SP}}(t,\rmu,\rmxi)=&g_0\left(1-\tau(a_1u+a_2v)+(t-\tau)A\right).
\end{split}
\end{equation}
Thus the computational cost can be significantly reduced and the accuracy is improved in smooth flow region as well.

Based on the above gas distribution functions Eq.(\ref{eq_f_dis}) and Eq.(\ref{eq_f_con}), a time-evolving flux function can be obtained according to Eq.(\ref{eq_f_F}). If the flux function is integrated within a time step directly, second-order time accuracy can be obtained with only one stage \textcolor{black}{since both the flux and its time derivative are involved}. However, if combined with the two-stage temporal discretization \cite{LPan2016II}, fourth-order time accuracy can be achieved in a highly efficient way, with the help of time derivative of the flux which can be simply computed from the gas-kinetic flux. This is introduced in the following.

\subsection{Fourth order Two-stage time integration}
Based on the unsteady gas-kinetic flux, an efficient two-stage time marching scheme can be constructed for the CPR framework to achieve fourth-order time accuracy by adopting the two-stage temporal discretization \cite{JLi2016, ZDu2018}. 
For the semi-discrete CPR method,  Eq.(\ref{eq_CPR}), the solution ${\rmQ}$ can be updated by, 
\begin{equation}\label{eq_S2O4}
\begin{split}
&{\rmQ}^*={\rmQ}^n+\frac12 \Delta t \mathcal{R}({\rmF}^n)+\frac18 \Delta t^2 \frac{\partial \mathcal{R}({\rmF}^n)}{\partial t},\\
&{\rmQ}^{n+1}={\rmQ}^n+\Delta t \mathcal{R}({\rmF}^n)+\frac16 \Delta t^2
\left(\frac{\partial \mathcal{R}({\rmF}^n)}{\partial t}+2\frac{\partial \mathcal{R}({\rmF}^*)}{\partial t}\right),
\end{split}
\end{equation}
where ${\rmF}^*={\rmF}({\rmQ}^*,t)$ is the flux at the middle stage $t^*=t^n+\Delta{t}/2$. It has been proved that the above two-stage temporal discretization achieves fourth-order time accuracy for hyperbolic conservation law. The success lies in the use of both the flux and its first-order time derivative. 

If the computational mesh does not change with time, $\mathcal{R}$ is a linear function of ${\rmF}$, thus $\partial \mathcal{R}({\rmF})/\partial t=\mathcal{R}(\partial{\rmF}/\partial t)$. So the key point is to compute the flux and its time derivative for both stages. First of all, the time-dependent flux can be obtained through the integration of the distribution function in the phase space according to Eq.(\ref{eq_f_F}), denoted as $\rmF(\rmQ^n,t)$. As it is not a linear function of $t$, a simple fitting method can be adopted to obtain the approximated flux and its first-order derivative. It is based on the time integration of $\rmF(\rmQ^n,t)$ within $[t_n,t_n+\delta]$,
\begin{equation}
\label{eq_S2O4_F(t)_int}
\widehat{\rmF}(\rmQ^n,\delta)=\int_{t^n}^{t^n+\delta}\rmF(\rmQ^n,t)\rmd{t}.
\end{equation}
As $\rmF(\rmQ^n,t)$ is approximated by a linear function $\tilde{\rmF}(\rmQ^n,t)=\rmF^n+(t-t^n)\partial_t\rmF^n$ within the time interval $[t^n,t^n+\Delta t]$, the time integration gives, 
\begin{equation}
\label{equ:S2O4_linear_F(t)_int_II}
\begin{split}
 &\frac12\Delta{t}\rmF^n+\frac18\Delta{t^2}\partial_t\rmF^n=\widehat{\rmF}(\rmQ^n,\Delta t/2),\\
 &\Delta{t}\rmF^n+\frac12\Delta{t^2}\partial_t\rmF^n=\widehat{\rmF}(\rmQ^n,\Delta t),
\end{split}
\end{equation}
where the left hand side is the time integration of $\tilde{\rmF}(\rmQ^n,t)$ within $[t^n,t^n+\Delta t/2]$ and $[t^n,t^n+\Delta t]$ respectively, while the right hand side is obtained according to Eq.(\ref{eq_S2O4_F(t)_int}). By solving the equation set Eq.(\ref{equ:S2O4_linear_F(t)_int_II}), we have
\begin{equation}
\label{equ:S2O4_linear_F(t)_int_III}
\begin{split}
 &\rmF^n=(4\widehat{\rmF}(\rmQ^n,\Delta t/2)-\widehat{\rmF}(\rmQ^n,\Delta t))/\Delta t,\\
 &\partial_t\rmF^n=4(\widehat{\rmF}(\rmQ^n,\Delta t)-2\widehat{\rmF}(\rmQ^n,\Delta t/2))/\Delta t^2.
\end{split}
\end{equation}  
Similarly, the approximated flux $\rmF^{\ast}$ and its time derivative $\partial_t\rmF^{\ast}$ can be obtained, by simply replacing the superscript $n$ in Eq.(\ref{equ:S2O4_linear_F(t)_int_III}) with $\ast$.

Now the final two-stage gas-kinetic CPR framework can be developed through Eqs.(\ref{eq_S2O4}) and (\ref{equ:S2O4_linear_F(t)_int_III}). At the first stage, with the solution at $t^n$, the gas distribution function is constructed at each SP and FP through Eqs.(\ref{eq_f_dis}) and (\ref{eq_f_con}), respectively. Then the time integrals $\widehat{\rmF}(\rmQ^n,\Delta t/2)$ and $\widehat{\rmF}(\rmQ^n,\Delta t)$ can be computed through Eq.(\ref{eq_S2O4_F(t)_int}), and thus the solution at each SP is updated by
\begin{equation}
\label{equ:S2O4_CPR_stage_I}
  \rmQ_{i,j}^*=\rmQ_{i,j}^n + \mathcal{R}_{i,j}(\widehat{\rmF}),
\end{equation}  
where the flux at each SP and FP is computed by
\begin{equation}
\label{equ:S2O4_CPR_stage_I_flux}
  \widehat{\rmF}=\widehat{\rmF}(\rmQ^n,\Delta t/2).
\end{equation}  

At the second stage, the time integration of the flux $\widehat{\rmF}(\rmQ^*,\Delta t/2)$ and $\widehat{\rmF}(\rmQ^*,\Delta t)$ at each SP and FP can be obtained in a similar way with the solution at $t^*$. Then the solution at each SP is updated by
\begin{equation}
\label{equ:S2O4_CPR_stage_II}
  \rmQ_{i,j}^{n+1}=\rmQ_{i,j}^* + \mathcal{R}_{i,j}(\widetilde{\rmF}),
\end{equation}   
where the flux at each SP and FP is computed by
\begin{equation}
\label{equ:S2O4_CPR_stage_II_flux}
  \widetilde{\rmF}=\frac{8}{3}\widehat{\rmF}(\rmQ^n,\Delta t/2) - \frac{1}{3}\widehat{\rmF}(\rmQ^n,\Delta t) -\frac{8}{3}\widehat{\rmF}(\rmQ^*,\Delta t/2)  + \frac{4}{3}\widehat{\rmF}(\rmQ^*,\Delta t) .
\end{equation}  

Compared with traditional CPR methods, which usually adopt the multi-stage R-K method, the current two-stage fourth-order gas-kinetic CPR can be more efficient because the number of stages is reduced from four or five to two, and the viscous flux computation is avoided.

\subsection{Hybrid CPR/SCFV approach for shock capturing}
 When simulating high-speed flows with shock waves, it is difficult for a high-order scheme with IDOFs to preserve its high resolution as the solution is continuous within the element. Numerical oscillations always occur due to the Gibbs phenomenon. Limiters or artificial viscosity can be adopted to suppress the oscillations. However, the subcell resolution may be lost because shock-capturing schemes are only first order $O({\delta} x)$  at the element containing the shock. 

To capture shock waves robustly and accurately, a hybrid CPR/SCFV approach is developed by extending the subcell finite volume (SCFV) limiting procedure \cite{MDumbser2014,MDumbser2016} to the CPR framework. The smooth flow regions are solved by the CPR method directly, while regions containing shock waves are solved by the SCFV method. With the hybrid  method, the subcell resolution can be preserved automatically and shock waves can be captured with high resolution and robustness. 

To construct a hybrid method, we need to identify the "troubled cells", and determine how the solutions are transferred between both methods. 
To distinguish regions containing shocks from smooth regions, a simple and effective shock detector is adopted to mark troubled cells \cite{GFu2017} at the beginning of each time step,
\begin{equation}
\label{equ:detector}
I_{\Omega_0}=\frac{\sum\limits_{l=1}^3 |\bar{Q}_0^n-\bar{Q}_{0,l}^n|}{\max\limits_{0\leq l \leq 3}\bar{Q}_l^n},
\end{equation} 
\textcolor{black}{where $\bar{Q}_l^n$ indicates the averaged density or total energy on the target cell with $l=0$, and its three face neighbors with $l=1\sim3$. $\bar{Q}_{0,l}^n$ is the average of polynomial $Q_l^n(\rmx)$ ($l=1\sim3$) on the target cell.}
\textcolor{black}{As recommened in Ref.\cite{GFu2017}, the target cell is marked as a troubled cell when $I_{\Omega_0}>0.12$.} However, in our study, it is found too low, and too many cells are marked in some cases. Thus, for all numerical tests below, the threshold is set as 0.3, and meanwhile the face neighbors of a troubled cell are also marked. The new strategy seems to work better. Nevertheless, we note that the current detector may not be the best. More investigations are necessary in the future. 

For the unmarked, the CPR framework is adopted directly to keep the high accuracy and efficiency in smooth flow regions. Meanwhile, a robust SCFV method is applied to those troubled cells, which need to be further partitioned into a set of subcells. For clarity, the original cell is also called a main cell. 

\begin{figure}[H]
  \centering
  \includegraphics[scale=0.25]{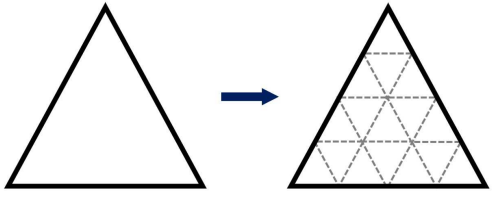}\\
  \captionsetup{labelsep=period}
  \caption{The partition of subcells.}
  \label{Subcells}
\end{figure}

Assume that $\Omega_i$ is marked as a troubled cell at $t_n$. For a fourth-order CPR with $k=3$, each cell contains 10 IDOFs. To preserve the subcell resolution and avoid the complexity of partition encountered in the SV method \cite{ZJWang2004}, $\Omega_i$ is uniformly partitioned into $(k+1)^2=16$ subcells $\{\Omega_{i,j}\},~j=1\sim16$, as shown in Fig.\ref{Subcells}.

The data transfer from CPR to SCFV is directly based on the solution polynomial. By taking average of $\rmQ_i^n(\rmx)$ on subcell $\Omega_{i,j}$, the subcell averaged solution can be obtained 
\begin{equation}
\label{eq_subcell_average}
\bar{\rmQ}_{i,j}^n=\frac1{|\Omega_{i,j}|}\int_{\Omega_{i,j}}\rmQ_i^n(\rmx)\rmd\Omega,~~~\forall\Omega_{i,j}\in\{\Omega_{i,j}\}.
\end{equation}
With the subcell averaged solutions, a robust second-order TVD finite volume reconstruction is implemented on each subcell $\Omega_{i,j}$ to obtain a non-oscillatory linear solution distribution, denoted as $\tilde{\rmQ}_{i,j}^n(\rmx)$. The reconstruction stencil for a subcell only involves its face-neighboring subcells. We note that, to provide the stencil of the reconstruction, the face-neighboring main cells of $\Omega_i$ also need to be partitioned into subcells, even if they are not marked as troubled cells. Details of the reconstruction can be referred to Ref.\cite{ZJWang2004,TJBarth1989}. As a result, the original cubic solution polynomial $\rmQ_i^n(\rmx)$ on the troubled cell $\Omega_i$ is replaced by a set of piecewise linear solution distributions, i.e., $\{\tilde{\rmQ}_{i,j}^n(\rmx)\}$.

With the new distributions, we can obtain the solutions at interfaces bounding each subcell, on which the flux is computed. Then the subcell averaged solution are updated in the finite volume framework. We note that the flux solver and the time-stepping method used on these subcells are consistent with that used in the CPR framework.

The data transformation from SCFV to CPR is actually a reconstruction. At $t^{n+1}$, with the set of subcell averaged solutions on $\Omega_i$, i.e., $\{\bar{\rmQ}_{i,j}^{n+1}\}$, the cubic solution polynomial $\rmQ_i^{n+1}(\rmx)$ can be recovered by a least-square reconstruction.
\begin{equation}
\label{equ_SCFV_t_n+1_LS}
\mathrm{min}\left[\sum_{j=1}^{16}\left(\frac{1}{|\Omega_{i,j}|} \int_{\Omega_{i,j}}\rmQ_i^{n+1}(\rmx)\rmd\Omega-\bar{\rmQ}_{i,j}^{n+1}\right)^2\right].
\end{equation}
To conserve the main cell averaged solution on $\Omega_i$ directly, the zero-mean basis function is adopted in $\rmQ_i^{n+1}(\rmx)$. Then the shock detector is applied to mark troubled cells at $t^{n+1}$. If $\Omega_i$ is marked again, $\{\bar{\rmQ}_{i,j}^{n+1}\}$ can be directly used for the limiting procedure at $t^{n+1}$. Otherwise, if $\Omega_i$ is not marked at $t^{n+1}$, the CPR framework is adopted on this cell. The solution at SPs can be computed with $\rmQ_i^{n+1}(\rmx)$. More details of the limiting procedure can be found in Ref.\cite{MDumbser2016}. 

Since each subcell is treated separately, discontinuities are introduced on the interfaces between subcells. As a result, shock waves can be resolved in the scale of subcells. The thickness of shock waves can be smaller than the size of main cells. On the contrary, for traditional limiters applied to main cells directly, discontinuities can only exist on interfaces between main cells, and shock waves can only be resolved in the scale of main cells. It is obvious that the current method is able to capture shock waves with higher resolution. In short, the hybrid CPR/SCFV approach fully combines the high accuracy and efficiency of CPR with the robustness and high resolution of SCFV.

\subsection{Remarks on efficiency}
In terms of the computational cost, the main difference between the present and the traditional CPR methods is the flux evaluation, and the time-stepping approach. To achieve fourth-order time accuracy, the traditional CPR usually adopts the widely used five-stage fourth-order R-K method. In contrast, the current scheme only needs two stages. 
Denote $T_{gks}$ as the computational cost in a time step for the present method, and $T_{cpr}$ for the traditional CPR. They can be expressed as 
\begin{equation}
\label{cost_gks}
\begin{split}
&T_{gks}=2(T_{gks,sp}^{inv+vis}+T_{gks,fp}^{inv+vis}+T_{gks}^{res}),\\
&T_{cpr}      =5(T_{cpr,sp}^{inv}+T_{cpr,sp}^{vis}+T_{cpr,fp}^{inv}+T_{cpr,fp}^{vis}+T_{cpr}^{res}),
\end{split}
\end{equation}
in which the factors 2 and 5 indicate the number of stages. The subscripts $sp$ and $fp$ represent the corresponding cost of flux computation at SPs and FPs. 
The superscript $res$ indicates the rest part of the computational cost, including reconstruction and residual evaluation. The superscripts $inv$ and $vis$ are for invscid and viscous parts, respectively. For the traditional CPR, the inviscid flux and viscous flux indicated by the superscripts $inv$ and $vis$ are computed separately. The flux at SPs is computed according to the flux function of the N-S equations. The inviscid flux at FPs is usually computed by Riemann solvers. For the gas-kinetic flux solver, these inviscid and viscous fluxes are coupled and computed simultaneously. 

Now we can easily estimate the computational cost of these two methods. For simplicity, only smooth flow is considered, thus no limiter is considered in the estimate. The invscid flux is computed with the Roe scheme in CPR and the viscous part is by the LDG scheme. Due to the higher order reconstruction for the viscous flux, $T_{cpr}^{res}$ is a little bit larger than $T_{gks}^{res}$ by about ten percent from numerical experiments. For the flux computation at FPs, the gas-kinetic solver is more expensive than the Roe scheme plus LDG, with a ratio $T_{gks,fp}^{inv+vis}/(T_{cpr,fp}^{inv}+T_{cpr,fp}^{vis})\approx 3.4$. For the flux at at SPs, the ratio between the continuous gas-kinetic solver and analytic N-S flux function becomes $T_{gks,sp}^{inv+vis}/(T_{cpr,sp}^{inv}+T_{cpr,sp}^{vis})\approx 2.1$. The ratios among different parts for GKS are nearly $T_{gks,sp}^{inv+vis}:T_{gks,fp}^{inv+vis}:T_{gks}^{res} \approx 1:2.5:1.1$. Therefore, it can be estimated that the overall computational cost between CPR and the current scheme is about $T_{cpr}/T_{gks}\approx 1.3$ when simulating viscous flows. Detailed data can be found in the following numerical test. 

Although the codes are not optimized, one can see that the efficiency improvement of the current method mainly comes from the adoption of the two-stage time stepping method, using only two stages rather than five stages in R-K. The straightforward computation of the time-derivative of the flux also contributed to the efficiency gain. Besides, it can also be observed that the overall cost of GKS is less than CPR, even if a four-stage R-K method is adopted. Furthermore, for supersonic flows, it is necessary to adopt a limiter for shock capturing, which usually takes a large mount of additional CPU time. As can be expected, the current method should can be more efficient. {\color{black}For the current hybrid method, a little additional CPU time is requrired for the data transfer between CPR and SCFV for troubled cells, but the adopted TVD limiter on subcells is low cost.} In addition, for three-dimensional flow simulations the ratio of $T^{res}$ also increases remarkably, thus using the two-stage time-stepping method can achieve higher efficiency.

\section{Numerical tests}\label{section3}
Several benchmark flows are simulated to validate the performance of the current scheme. The ratio of specific heats is $\gamma=1.4$ in all of these tests. The collision time for viscous flows is computed by
\begin{equation}\label{eq_tau}
  \tau =\frac{\mu}{p}, \quad
  \tau_n =\tau+\epsilon_2\left|\frac{p^L-p^R}{p^L+p^R}\right|^{\epsilon_3}\Delta t.
\end{equation}
The variable $\tau_n$ is used to replace the physical collision time $\tau$ in the exponential function in Eq.(\ref{eq_BGK_exact}), for better controlling the numerical dissipation through the transition from $f_0$ to $g$, such as the second term in $\tau_n$. Here $p^L$ and $p^R$ are the pressure at the left and right sides of a cell interface. The coefficient $\epsilon_2=10$ is chosen. Another coefficient $\epsilon_3=\min(1,5\rma)$ is included to improve the accuracy when solving flows with very low Mach number. For inviscid flows, the physical viscosity $\mu/p$ in $\tau$ is replaced with $\epsilon_1\Delta t$, where $\epsilon_1=0.005$ is adopted in the current study.

To compute the time step $\Delta t$ based on the main cell size $h$, the CFL number is fixed to 0.1 in all tests. Boundary conditions are implemented with the help of ghost cells. When evaluating the computational cost, the simulation is carried out on the Intel Core i7-3770 CPU $\&$ 3.40 GHz. For convenience, the current scheme is denoted as CPR-GKS in the following.

For a high-order accurate scheme on triangular meshes, two questions are often encountered when presenting the numerical results. The first one is related to the IDOFs. For a numerical scheme with IDOFs, the variables at these IDOFs are meaningful which is different to a FV method.
In the present study, the errors in the accuracy test are computed directly from the variables at SPs,
\begin{equation}\label{eq_L1L2error}
\begin{split}
  L_2~\mathrm{error}=\sqrt{\frac{\sum_{i=1}^{N}\sum_{j=1}^{10}(q_{i,j}-q_{i,j}^e)^2}{10N}},
\end{split}
\end{equation}
in which $q_{i,j}$ and $q_{i,j}^e$ denote the numerical and analytical solution, respectively. $N$ is the number of cells.  Similarly, the residuals for steady flow problems are calculated, in which the steady state is achieved when the residual of velocity is less than $10^{-14}$ for the Couette flow, while $10^{-9}$ for the lid-driven cavity flow. Besides, all 2D contours are presented based on the averaged solution on subcells. 

\begin{figure}[H]
  \centering
  \begin{varwidth}[t]{\textwidth}
  \vspace{0pt}
  \includegraphics[scale=0.23]{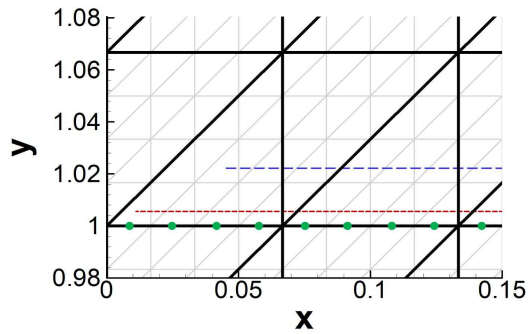}\\
  \end{varwidth}
  \captionsetup{labelsep=period}
  \caption{Enlarged view of the mesh with $h=1/15$  near the centerline $y=1$ in Shu-Osher problem.}
  \label{fig:Shu_Osher_Mesh}
\end{figure}

The second question is how to obtain the variables at given locations, especially along a specific line in 1D flows. Usually the interpolation can not be avoided which may lead to a result different to the original one. Here for simplicity, the computational triangular mesh cells are obtained from the triangulation of rectangular cells, except for the double Mach reflection flow and the viscous shock tube flow. Thus the interpolation can be simpler. For example, as shown in Fig.\ref{fig:Shu_Osher_Mesh}, where the black thick lines indicate main cells and the grey thin lines subcells, the variables at the equidistant green points can be computed directly by the solution polynomial Eq.(\ref{eq_Q(x,y)}) at both sides. Then the unique value can be determined through a simple arithmetic average of both sides. In the Shu-Osher flow, the averaged variables of a main cell and the corresponding subcells are also discussed, where the cell centroids can be easily found along the red dashed line and the blue one, respectively. 
Furthermore, based on the rectangular mesh cell the scale of cells can be well controlled which is important for accuracy tests. 

\subsection{Compressible Couette flow}
To quantitatively validate the accuracy of the current scheme, the compressible Couette flow is simulated, which has an analytical solution \cite{KXu2001,CZhang2018}. It is a steady flow between two parallel plates, driven by the upper plate, which has a constant speed $U_1=0.5$ and temperature $T_1=1$. The lower plate is stationary and adiabatic. The computational domain is $[0,4H]\times[0,2H]$ with $H=1$. The computational mesh is obtained by the triangulation of a rectangular mesh.
The Mach number is set as $\mathrm{Ma}=U_1/\sqrt{\gamma R T_1}=0.5$ and the Reynolds number is $\mathrm{Re}=\rho_1 U_1 H/\mu_1=500$ with $\rho_1=1$. The viscosity of the flow is determined by the linear law $\mu=\mu_1 T/T_1$. The Prandtl number is $\mathrm{Pr}=1$.
The flow variables at ghost cells are fixed to the analytical solutions for all boundaries.

For comparison, the results obtained by the traditional CPR method (denoted as CPR-LDG) is also presented, in which the inviscid flux is computed by the Roe scheme, while the viscous flux is computed by the LDG scheme, and the five-stage fourth-order R-K method is used for time stepping. The errors and convergence orders of density are presented in Table~\ref{AccuracyCouetteRho}. Both of the two schemes achieve the designed order of accuracy. There is only slight difference between them in terms of both errors and accuracy orders. An additional flow with very low Mach number $\mathrm{Ma}=0.02$ is also simulated to verify the performance of the current scheme in nearly incompressible flows. The results are shown in Table~\ref{AccuracyCouetteRho_II} in which the finest grid is not considered as the error is so low that it is affected by the round-off error. The designed order of accuracy can be achieved as well, and the errors significantly decrease with the Mach number. {\color{black}It should be noted that the fourth-order temporal accuracy of the current scheme is also validated through the inviscid isentropic vortex propagation flow, which is not presented here for brevity.} 

\begin{table}[H]
\centering
\captionsetup{labelsep=period}
\caption{Accuracy test in compressible Couette flow with $\mathrm{Ma}=0.5$}
\label{AccuracyCouetteRho}
\renewcommand\arraystretch{1.3}
\begin{tabular}{lllll}
\hline
$h$    & \multicolumn{2}{l}{CPR-GKS} & \multicolumn{2}{l}{CPR-LDG} \\ \cline{2-5}
     & $L_2$ error    & Order      & $L_2$ error    & Order      \\ \hline
1.0  & 3.24E-06       &            & 2.55E-06       &            \\
0.5  & 1.66E-07       & 4.29       & 1.53E-07       & 4.06       \\
0.25  & 9.90E-09       & 4.07       & 9.52E-09       & 4.01       \\
0.125 & 6.01E-10       & 4.04       & 5.98E-10       & 3.99      \\ \hline
\end{tabular}
\end{table}

\begin{table}[H]
\centering
\captionsetup{labelsep=period}
\caption{Accuracy test in compressible Couette flow with $\mathrm{Ma}=0.02$}
\label{AccuracyCouetteRho_II}
\renewcommand\arraystretch{1.3}
\begin{tabular}{lllll}
\hline
$h$    & \multicolumn{2}{l}{CPR-GKS} & \multicolumn{2}{l}{CPR-LDG} \\ \cline{2-5}
     & $L_2$ error    & Order      & $L_2$ error    & Order      \\ \hline
1.0  & 5.76E-12       &            & 5.43E-12       &            \\
0.5  & 3.93E-13       & 3.87       & 3.23E-13       & 4.07        \\
0.25  & 2.55E-14       & 3.95       & 2.08E-14       & 3.95     \\ \hline
\end{tabular}
\end{table}

For a comparison of the efficiency, Table~\ref{CPU_Couette} shows the computational cost for 1000 time steps with $h=0.125$. We note that, for CPR-LDG, the cost of flux is the sum of inviscid flux and viscous flux which are presented separately. For example, the flux computation at SPs takes 3.29s for inviscid flux and 7.68s for viscous flux respectively. Thanks to the simplification of the gas distribution function at SPs, CPR-GKS takes less cost than CPR-LDG for the flux at SPs. For the flux at FPs, CPR-GKS is less efficient than CPR-LDG since the gas distribution function is more complicated. In consistent with the foregoing discussion, the overall cost for flux computation is close to each other. However, CPR-LDG takes much more cost than CPR-GKS in other parts which roughly includes the cost of the reconstruction, computing the flux divergence and flux correction terms. This is mainly due to the use of different time-stepping methods. With the efficient two-stage time-stepping method, CPR-GKS saves the cost of these parts in extra stages needed in CPR-LDG. As a result, CPR-GKS is more efficient than CPR-LDG. Furthermore, when solving high speed flows with shock waves, the limiter also takes a large part of the computational cost. Then it can be expected that the current scheme can be more efficient than traditional CPR, thanks to the two-stage flux evolution. The results of this case demonstrate the high accuracy and efficiency of the current scheme in compressible viscous flows. 

\begin{table}[H]
\centering
\captionsetup{labelsep=period}
\caption{Computational cost for 1000 time steps with $h=0.125$.}
\label{CPU_Couette}
\renewcommand\arraystretch{1.3}
\begin{tabular}{lllll}
\hline
Scheme  & Flux at SPs & Flux at FPs & Other parts & Total cost \\ \hline
CPR-GKS & 9.41 s      & 23.46 s      & 10.05 s       & 42.92 s     \\
CPR-LDG & 3.29 s + 7.68 s   & 5.81 s + 11.46 s      & 27.66 s      & 55.90 s     \\
Ratio   & 1.17        & 0.74        & 2.75        & 1.30       \\ \hline
\end{tabular}
\end{table}

\subsection{Shu-Osher Problem}
The Shu-Osher problem \cite{QWShu1989} involves the interaction between a shock wave and an entropy wave, which is simulated to validate the high resolution of the current scheme. The computational domain is $[0,10]\times[0,2]$ and two mesh sizes are considered with $h=1/15$ and $h=1/30$, respectively. 
The initial condition is
\begin{equation}\label{eq_Shu_Osher_initial}
(\rho,U,V,p)=\begin{cases}
(3.857134,2.629369,0,10.33333),& 0\leq x \leq 1, \\
(1+0.2\sin(5x),0,0,1),& 1< x \leq 10.
\end{cases}
\end{equation}

\begin{figure}[H]
  \centering
  \begin{varwidth}[t]{\textwidth}
  \vspace{0pt}
  \includegraphics[scale=0.33]{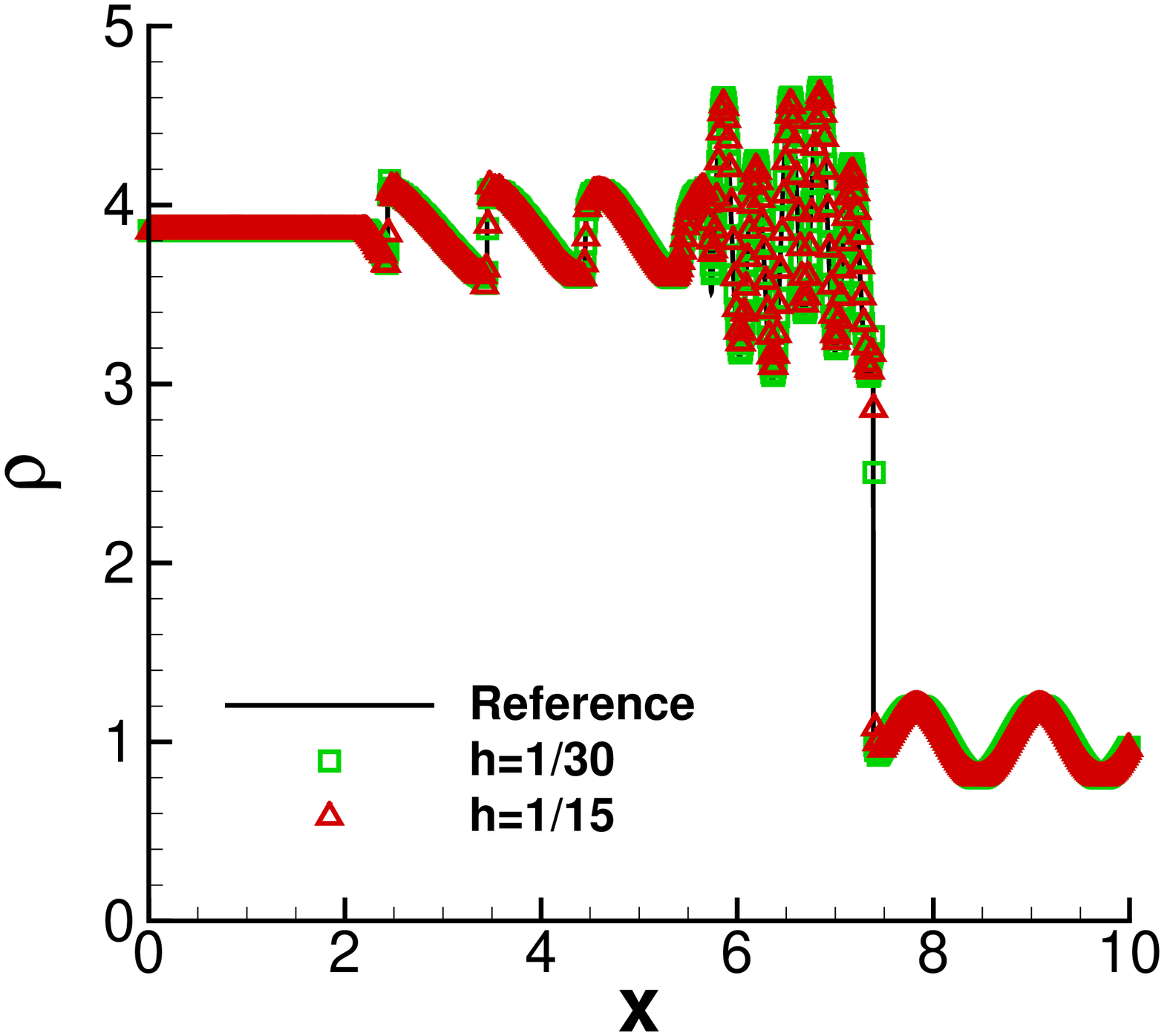}\\
  \end{varwidth}
  \qquad
  \begin{varwidth}[t]{\textwidth}
  \vspace{0pt}
  \includegraphics[scale=0.33]{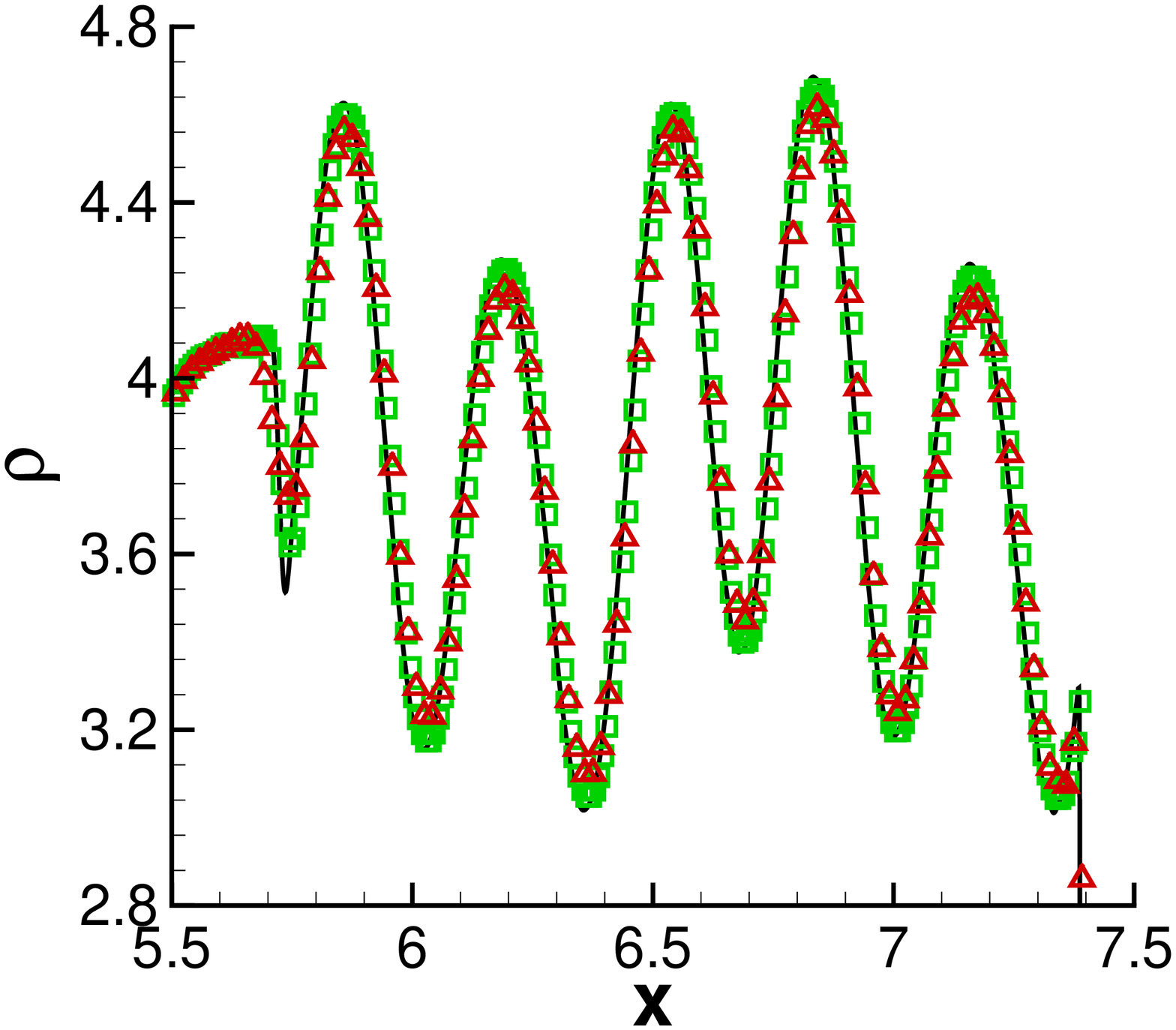}\\
  \end{varwidth}
  \captionsetup{labelsep=period}
  \caption{Density distribution along the horizontal centerline at t=1.8 in Shu-Osher problem.}
  \label{fig:Shu_Osher_rho_1D}
\end{figure}

\begin{figure}[H]
  \centering
  \begin{varwidth}[t]{\textwidth}
  \vspace{0pt}
  \includegraphics[scale=0.35]{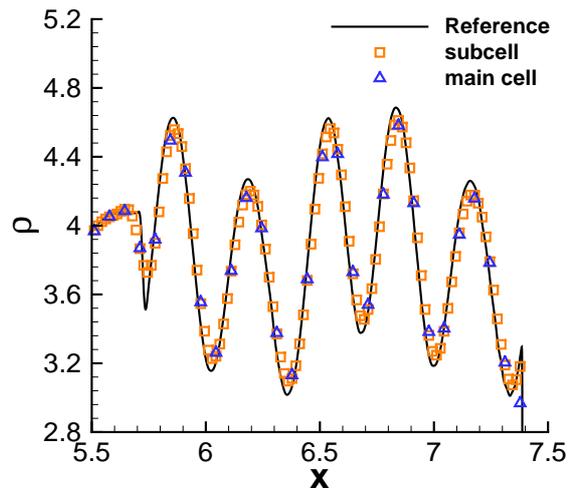}\\
  \end{varwidth}
  \captionsetup{labelsep=period}
  \caption{Cell-averaged density distribution of both subcells and main cells with $h=1/15$ at t=1.8 in Shu-Osher problem.}
  \label{fig:Shu_Osher_1D_Uave}
\end{figure}

\begin{figure}[H]
  \centering
  \begin{varwidth}[t]{\textwidth}
  \vspace{0pt}
  \includegraphics[scale=0.45]{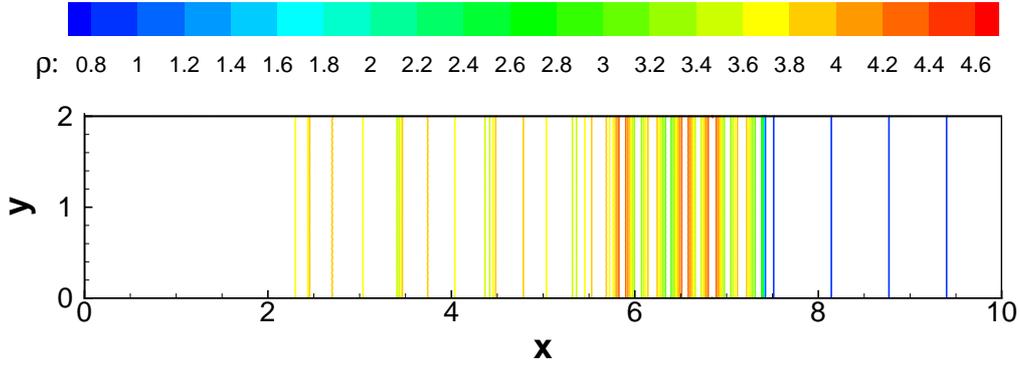}\\
  \end{varwidth}
  \captionsetup{labelsep=period}
  \caption{Density contours at t=1.8 with $h=1/15$ in Shu-Osher problem. 20 contours are drawn from 0.8 to 4.6.}
  \label{fig:Shu_Osher_contours_2D}
\end{figure}

\begin{figure}[H]
  \centering
  \begin{varwidth}[t]{\textwidth}
  \vspace{0pt}
  \includegraphics[scale=0.35]{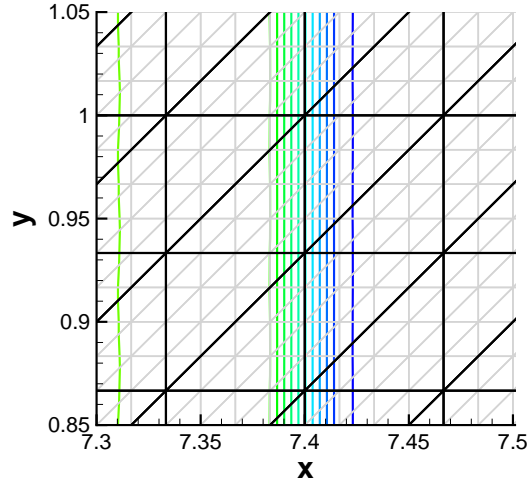}\\
  \end{varwidth}
  \captionsetup{labelsep=period}
  \caption{Density contours near the shock wave at $t=1.8$ with $h=1/15$ in Shu-Osher problem.}
  \label{fig:Shu_Osher_thickness}
\end{figure}

Fig.\ref{fig:Shu_Osher_rho_1D} presents the density distribution at $t=1.8$ along the horizontal centerline $y=1$. The results computed by the 1D HGKS \cite{QBLi2010} with $h=1/1000$ is chosen for reference. It can be observed that the current scheme can capture the density profile smoothly without numerical oscillation. With the mesh size $h=1/30$, the result obtained by CPR-GKS matches well with the reference data.
For a numerical scheme with IDOFs, the variables at these IDOFs are meaningful which is different to a FV method.  Here the cell averaged densities for $h=1/15$ are also presented for comparison (see Fig.\ref{fig:Shu_Osher_1D_Uave}). It can be observed that the main-cell averaged solutions are as accurate as the subcell counterpart. However, the subcell averaged solutions can provide much more details of the flow distribution with higher resolution.  Fig.\ref{fig:Shu_Osher_contours_2D} shows the 2D density contours. In particular, an enlarged view near the shock wave is also presented in Fig.\ref{fig:Shu_Osher_thickness}. The thickness of the predicted shock wave is less than the size of main cells, demonstrating the subcell resolution of the current scheme.

\begin{figure}[H]
  \centering
  \begin{varwidth}[t]{\textwidth}
  \vspace{0pt}
  \includegraphics[scale=0.5]{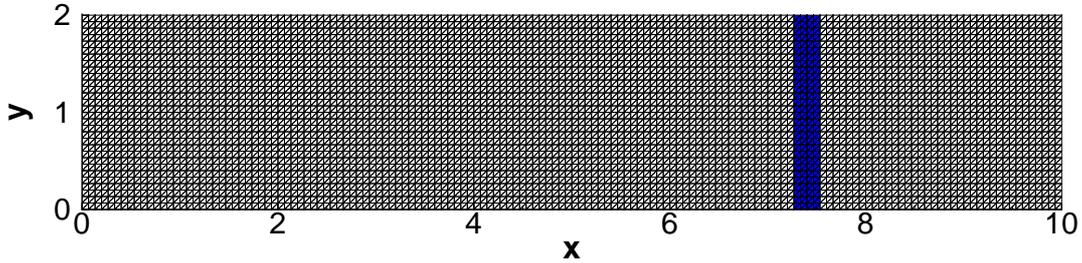}\\
  \end{varwidth}
  \captionsetup{labelsep=period}
  \caption{Distribution of troubled cells at $t=1.8$ with $h=1/15$ in Shu-Osher problem.}
  \label{fig:Shu_Osher_TBC}
\end{figure}

Besides, the troubled cells at $t=1.8$ are shown in Fig.\ref{fig:Shu_Osher_TBC}. It can be observed that the shock detector successfully marks the shock wave region which is simulated by the SCFV method. Nevertheless, it is noted that the current shock detecting strategy may not be the best. Furthermore, on these subcells the second-order TVD limiter is applied which usually introduces large numerical dissipation and may impair the resolution of wave structure when the shock wave is interacting with the entropy wave. A high-order limiter may be adopted to to improve the accuracy in the near-shock region and reduce the dependence of shock detectors.

\subsection{Blast wave problem}
\begin{figure}[H]
  \centering
  \begin{varwidth}[t]{\textwidth}
  \vspace{0pt}
  \includegraphics[scale=0.33]{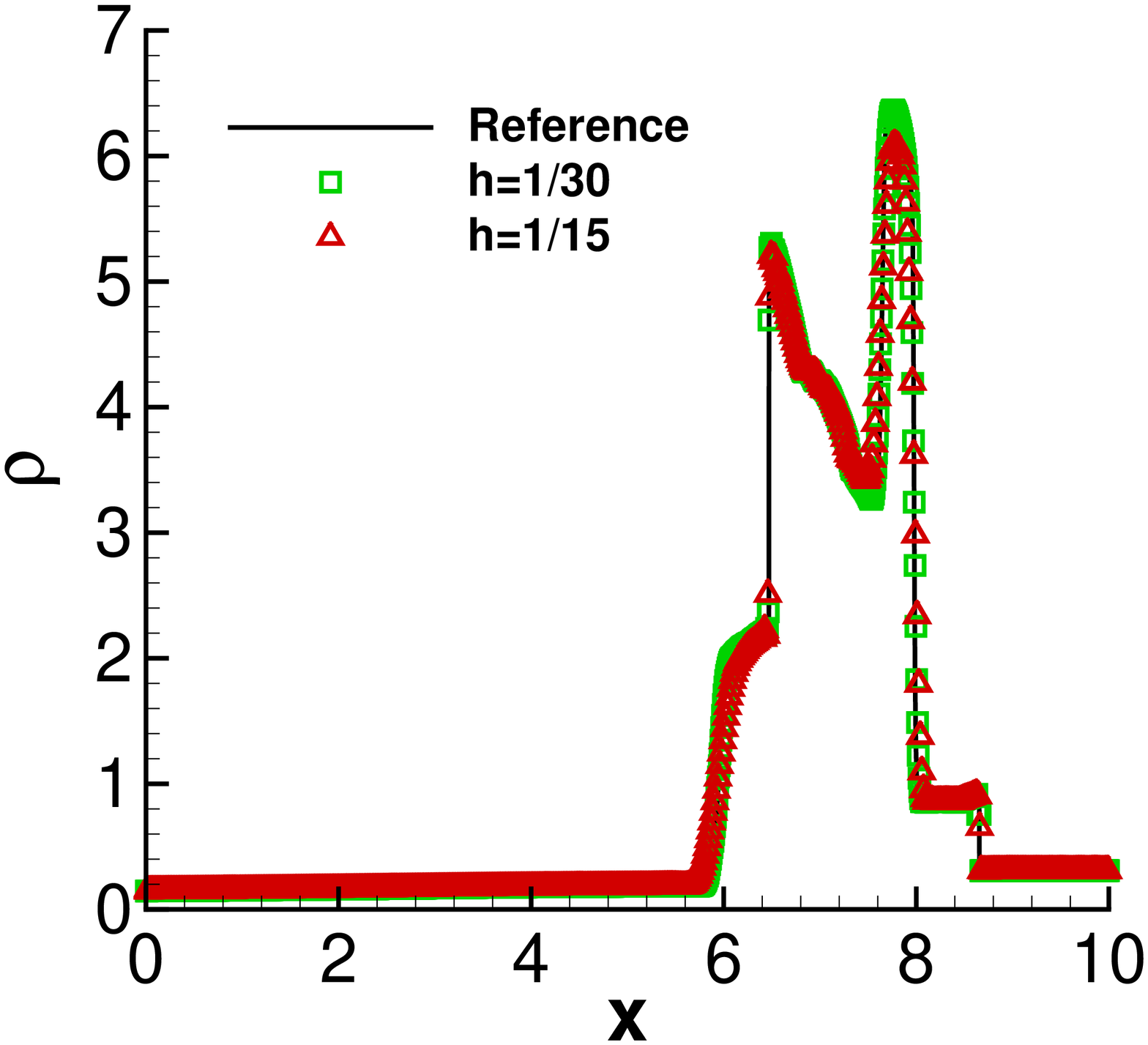}\\
  \end{varwidth}
  \qquad
  \begin{varwidth}[t]{\textwidth}
  \vspace{0pt}
  \includegraphics[scale=0.33]{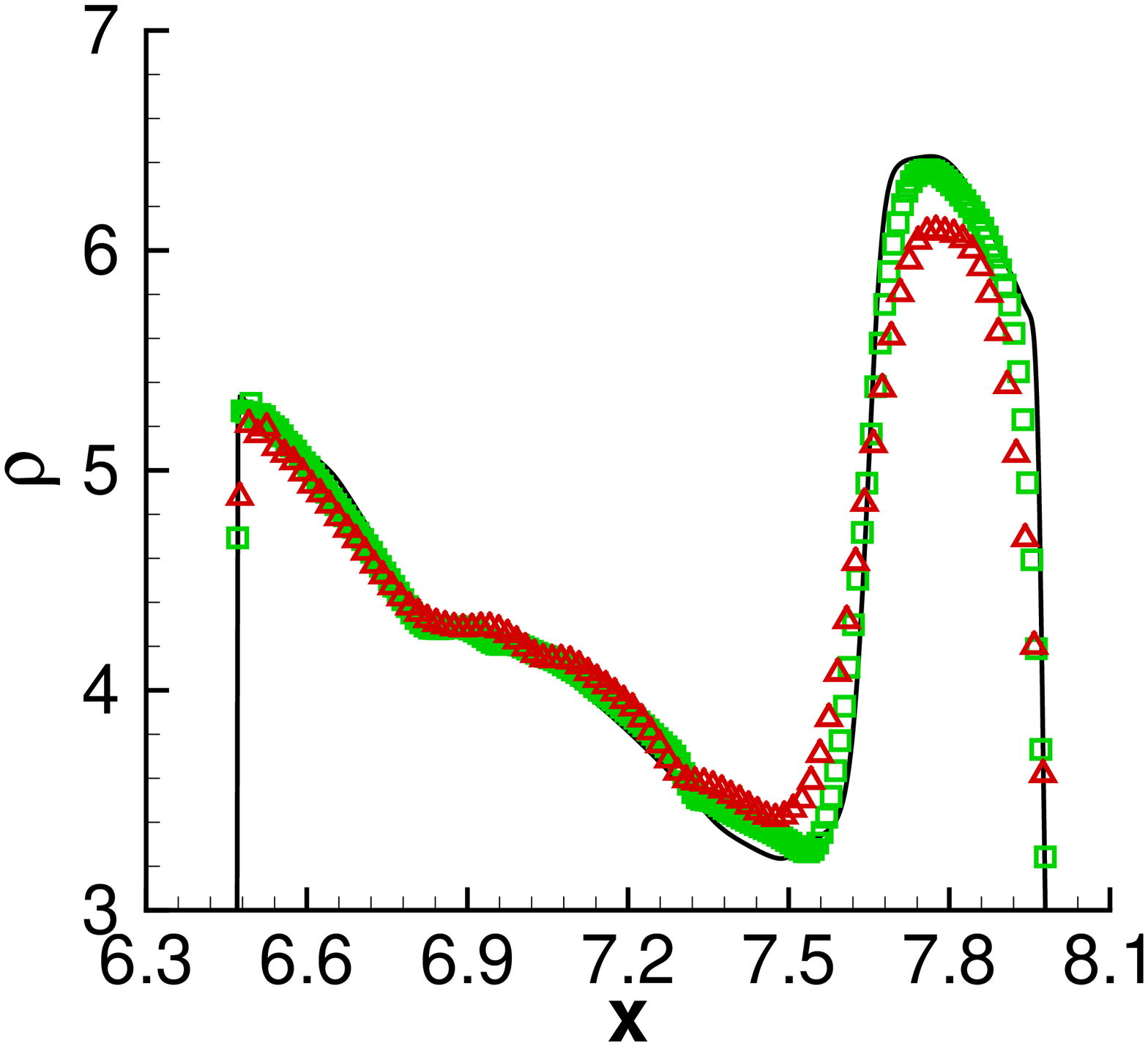}\\
  \end{varwidth}
  \captionsetup{labelsep=period}
  \caption{Density distribution along the horizontal centerline at t=0.38 in blast wave problem.}
  \label{fig:blast_rho_1D}
\end{figure}

The blast wave problem \cite{PWoodward1984} is a typical case with extremely strong shock waves, which can be used to validate the robustness and resolution of strong discontinuity. The computational domain is $[0,10]\times[0,2]$. The computational mesh is the same as that used in the Shu-Osher problem. The initial condition is
\begin{equation}\label{eq_Blast_initial}
(\rho,U,V,p)=\begin{cases}
(1,0,0,1000),& 0\leq x \leq 1, \\
(1,0,0,0.01),& 1< x \leq 9, \\
(1,0,0,100),& 9< x \leq 10.
\end{cases}
\end{equation}
The results computed by the 1D HGKS \cite{QBLi2010} with the mesh size $h=1/1000$ are chosen as the reference data. The density distribution at $t=0.38$ along the horizontal centerline $y=1$ is shown in Fig.\ref{fig:blast_rho_1D}. The current scheme successfully captures the strong shock wave interaction without oscillation. With the mesh size $h=1/30$, the density distribution matches very well with the reference data. Due to the extremely strong shock wave, it is very difficult to suppress the oscillation by applying limiters to the high-order solution polynomial directly. Although the results can be non-oscillatory as well, too much numerical dissipation is introduced leading to a much more dissipative flow distribution even with a fine mesh. In contrast, with the SCFV limiting procedure, the oscillation can be suppressed successfully. Moreover, since discontinuity can exist between subcells, shock waves are captured much more sharply. Compared to a traditional CPR using the MLP limiter and with the same mesh size $h=1/30$ \cite{JSPark2016}, the density profile obtained by the current scheme is much better. These results demonstrate the robustness and high resolution of the current scheme for strong shock waves.

\subsection{Double Mach reflection}
\begin{figure}[H]
  \centering
  \begin{varwidth}[t]{\textwidth}
  \vspace{0pt}
  \includegraphics[scale=0.28]{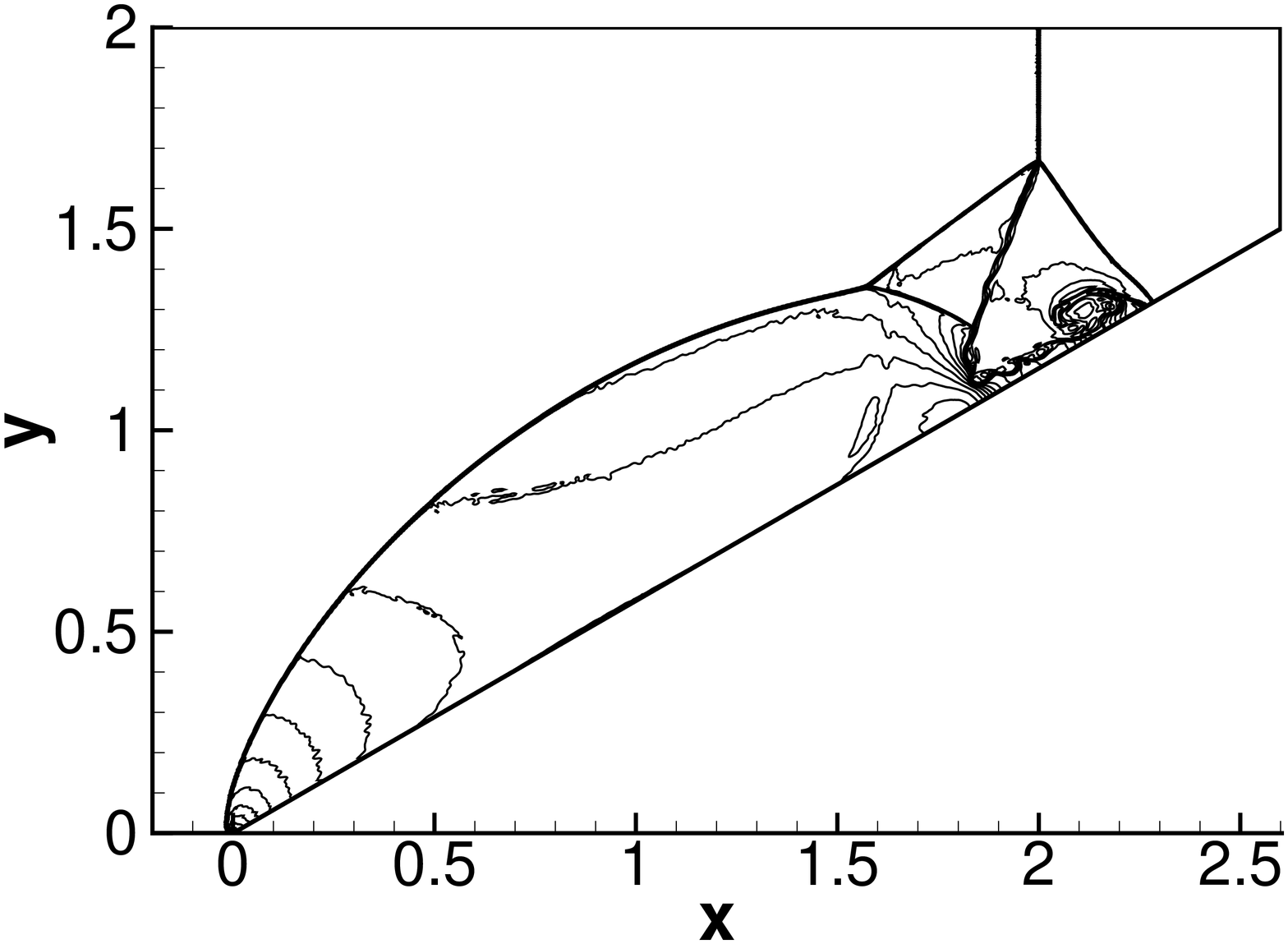}\\
  \end{varwidth}
  \qquad
  \begin{varwidth}[t]{\textwidth}
  \vspace{0pt}
  \includegraphics[scale=0.28]{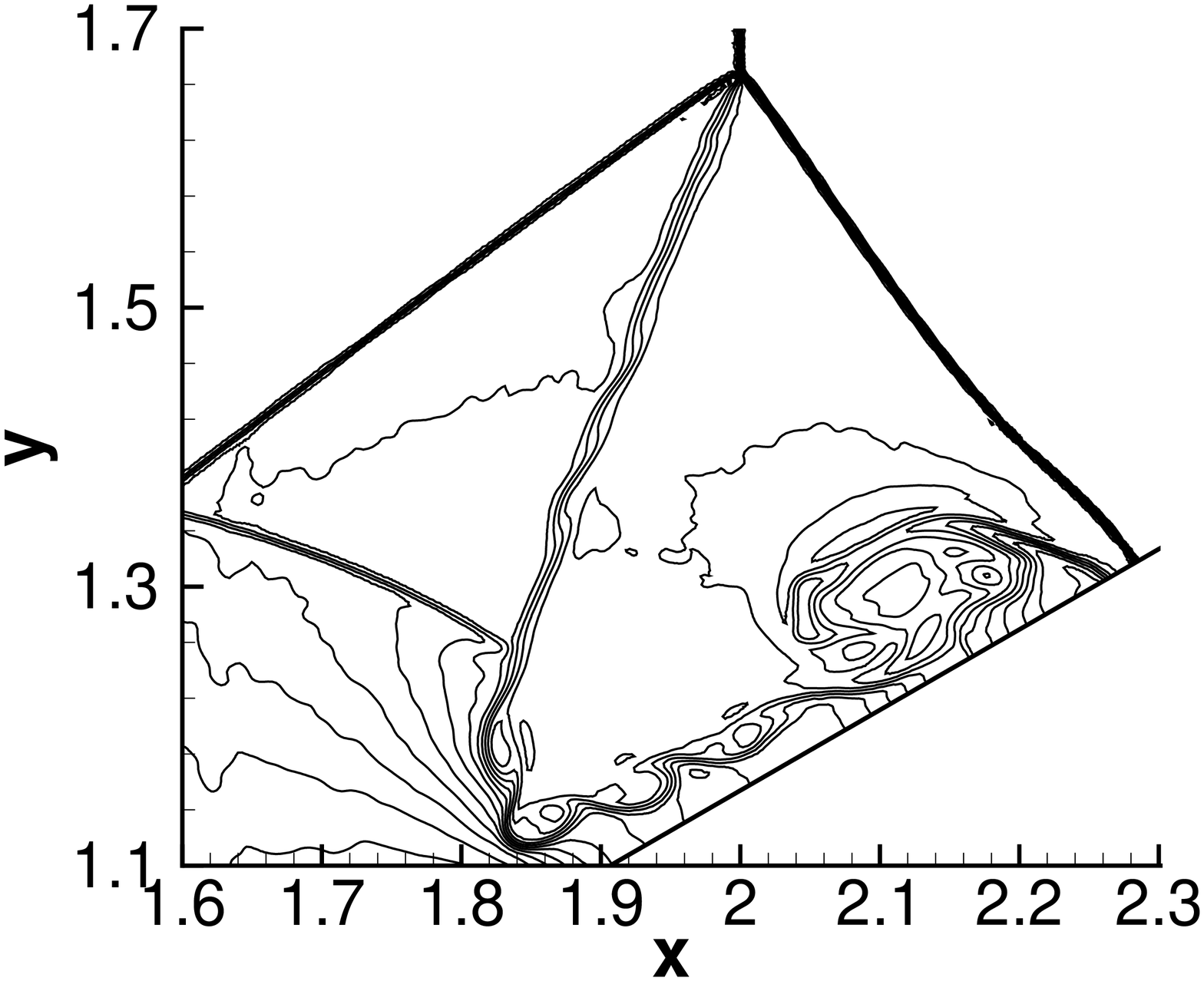}\\
  \end{varwidth}
  \begin{varwidth}[t]{\textwidth}
  \vspace{0pt}
  \includegraphics[scale=0.28]{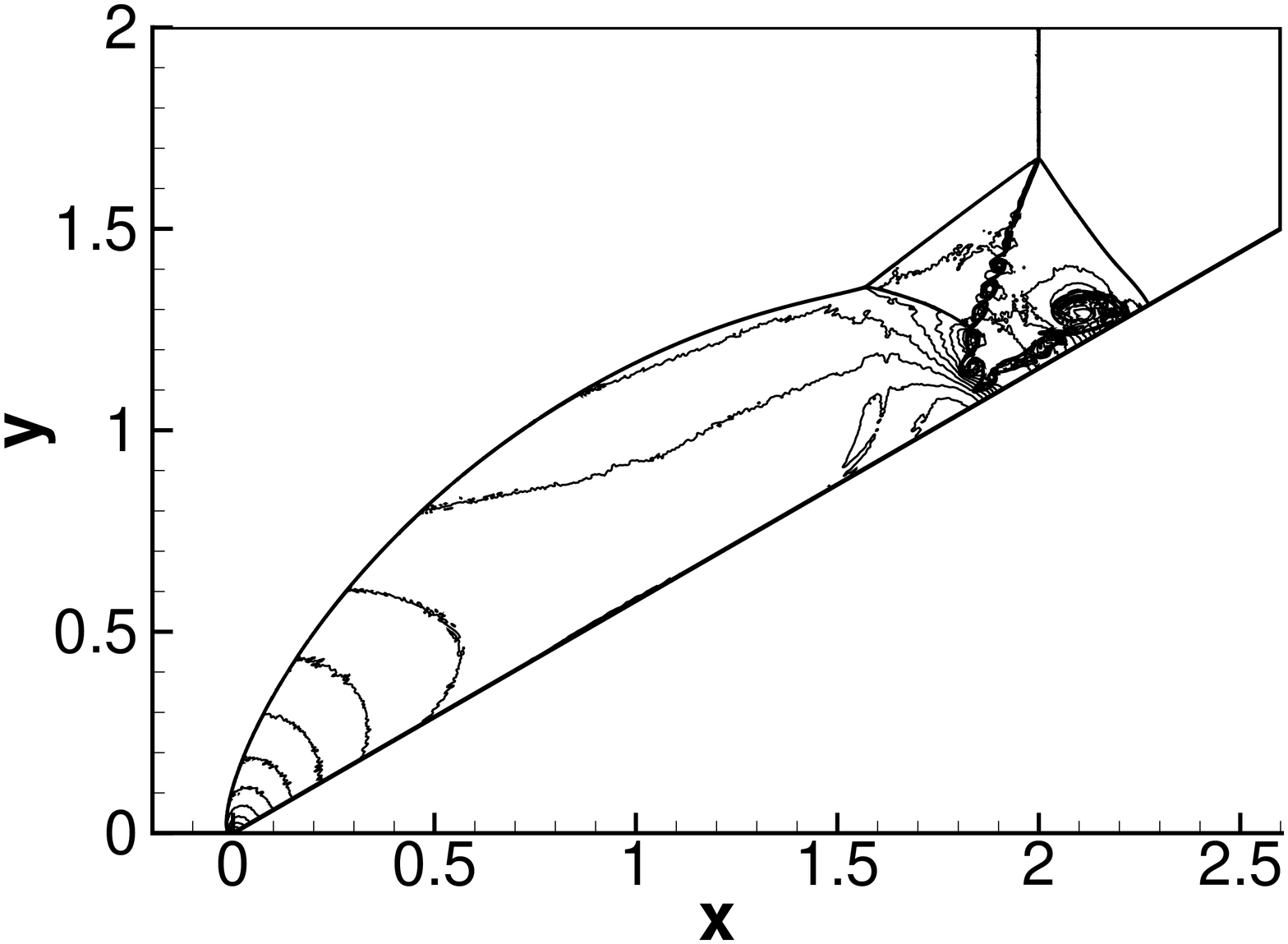}\\
  \end{varwidth}
  \qquad
  \begin{varwidth}[t]{\textwidth}
  \vspace{0pt}
  \includegraphics[scale=0.28]{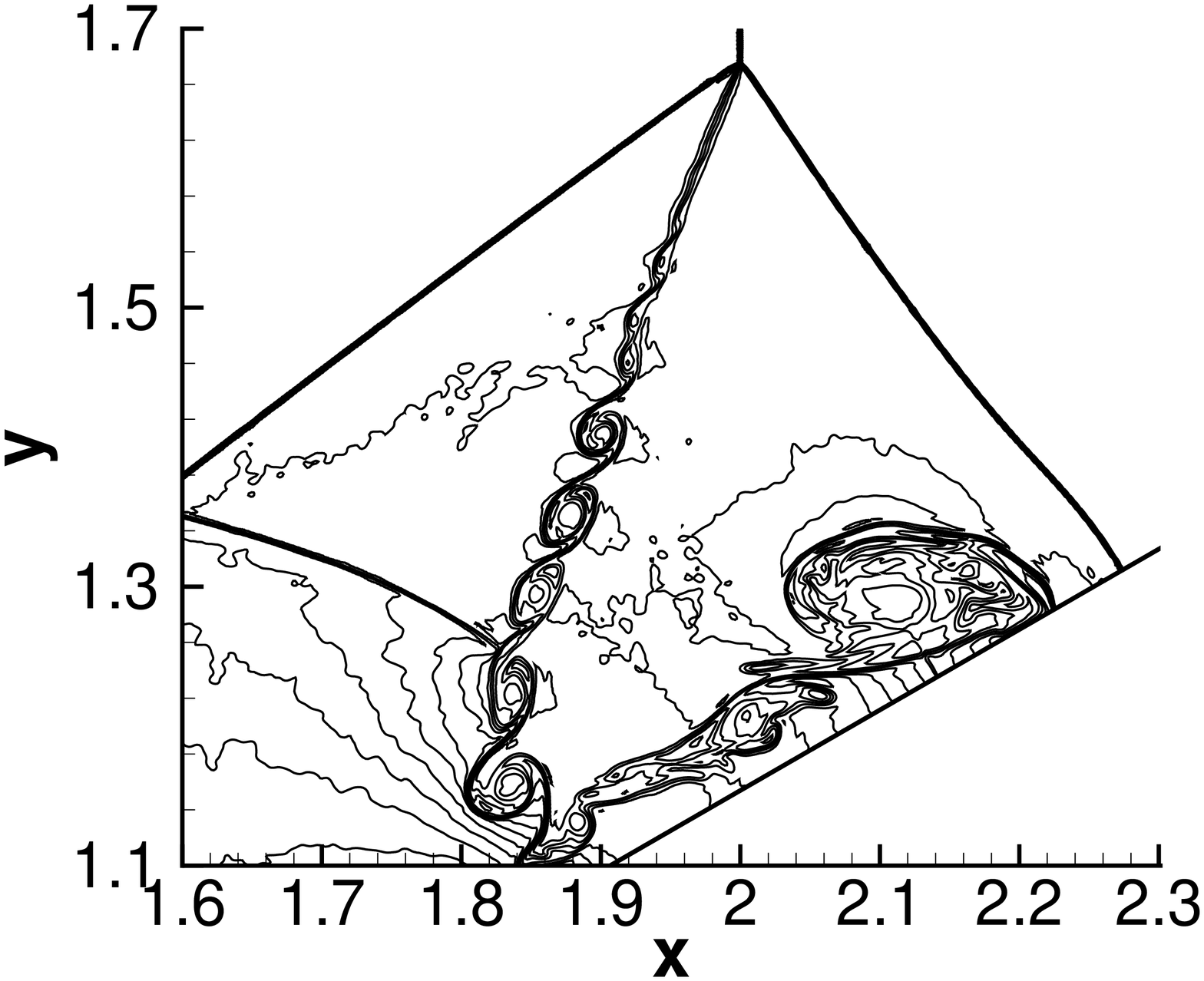}\\
  \end{varwidth}
  \captionsetup{labelsep=period}
  \caption{Density distributions at $t=0.2$ and the enlarged view near the Mach stem in double Mach reflection with mesh size $h=1/80$ (top) and $h=1/160$ (bottom). 30 contours are drawn from 2.0 to 22.5.}
  \label{fig:DMR_rho_contours}
\end{figure}

\begin{figure}[H]
  \centering
  \begin{varwidth}[t]{\textwidth}
  \vspace{0pt}
  \includegraphics[scale=0.32]{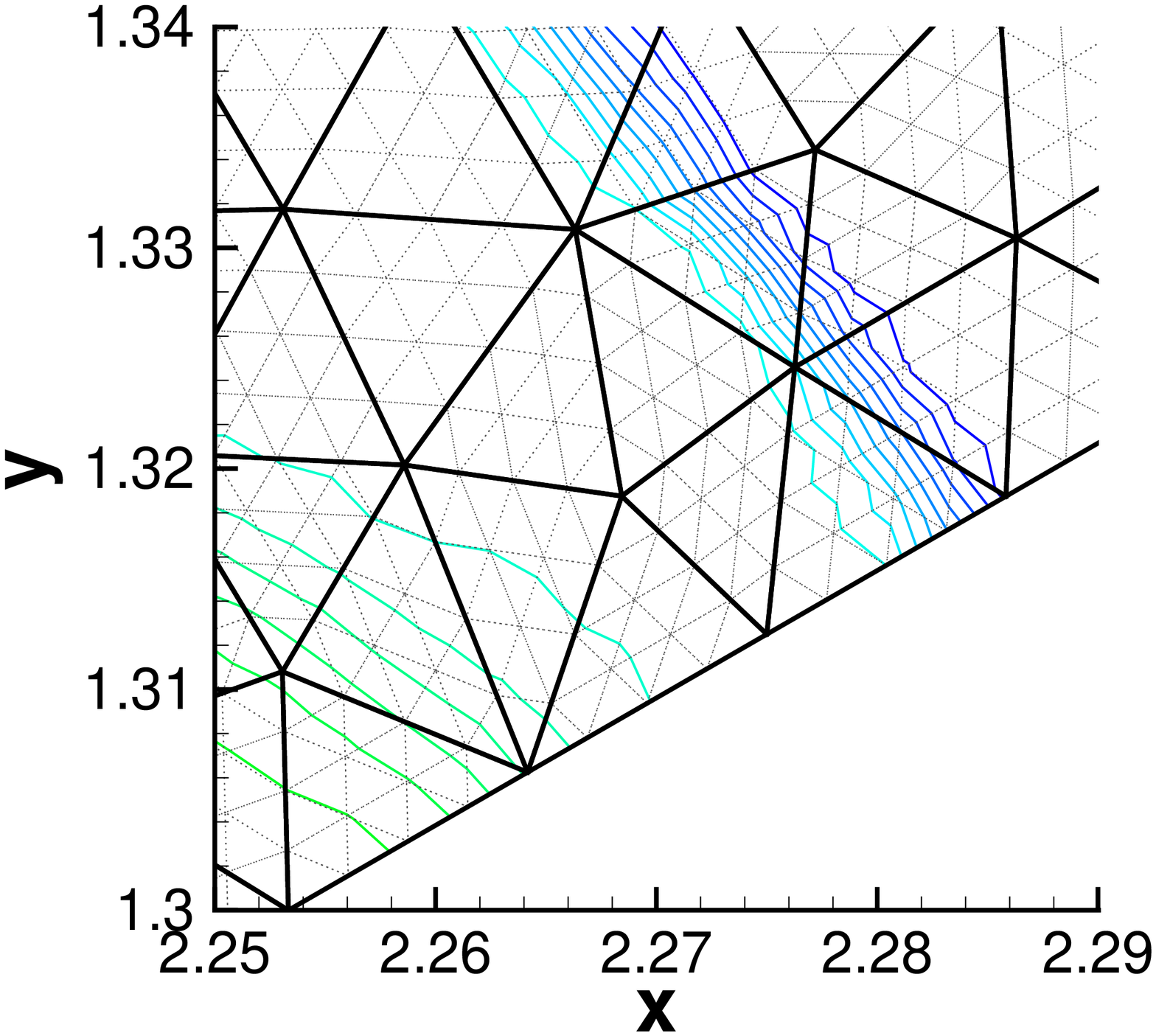}\\
  \end{varwidth}
  \qquad
  \begin{varwidth}[t]{\textwidth}
  \vspace{0pt}
  \includegraphics[scale=0.32]{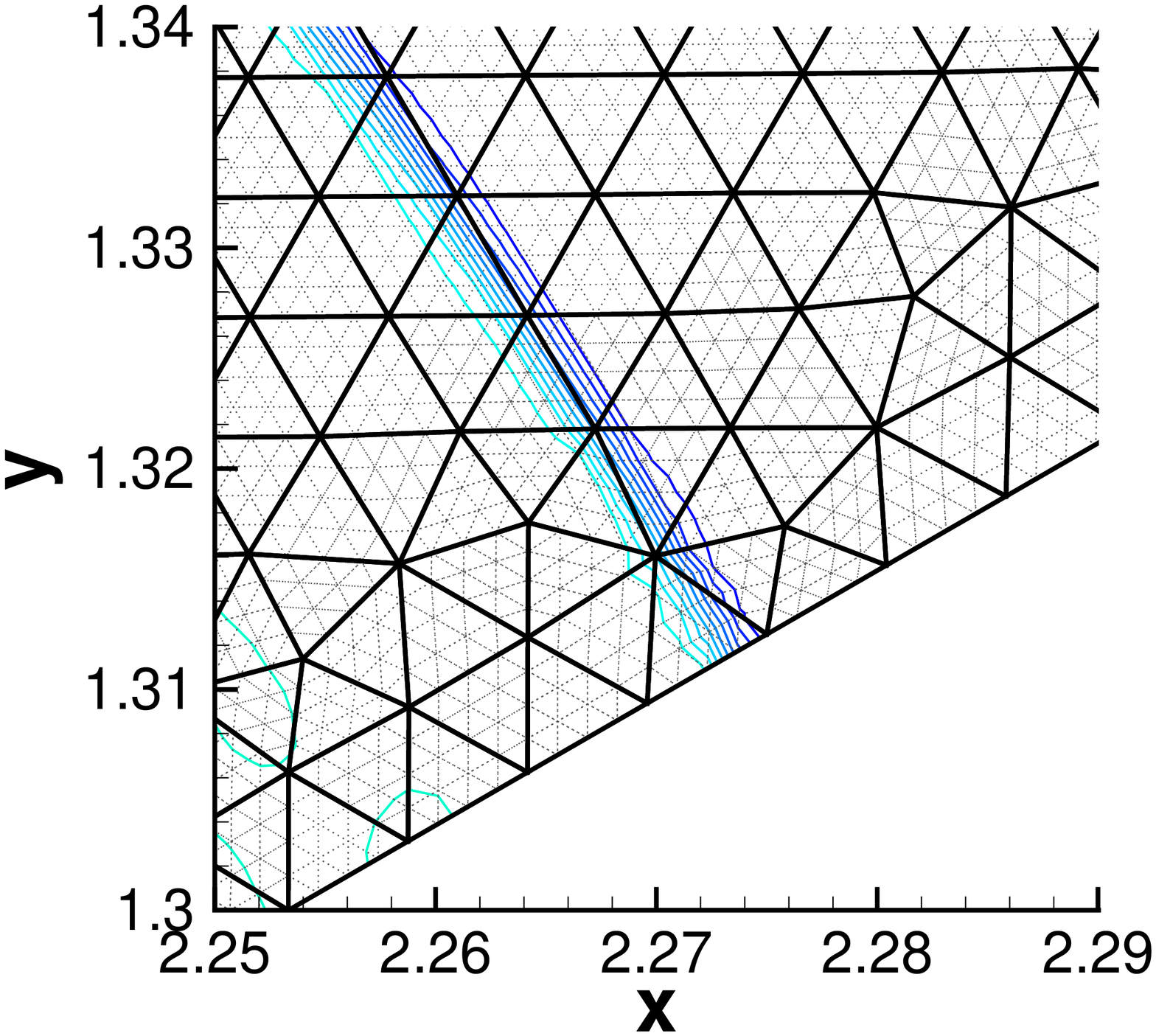}\\
  \end{varwidth}
  \captionsetup{labelsep=period}
  \caption{Density contours of the Mach stem near the wall with mesh size $h=1/80$ (left) and $h=1/160$ (right).}
  \label{fig:DMR_shock_thickness}
\end{figure}

The double Mach reflection \cite{PWoodward1984} is a 2D benchmark flow to validate the robustness and resolution of a numerical scheme in hypersonic flows. Initially, a right-moving shock wave with the Mach number $\mathrm{Ma}=10$ is located at $x=0$ with the condition
\begin{equation}\label{eq_DMR_initial}
(\rho,U,V,p)=\begin{cases}
(1.4,0,0,1),& x \leq 0, \\
(8,8.25,0,116.5),& x> 0,
\end{cases}
\end{equation}
which impinges on a $30^{\circ}$ wedge and leads to the double Mach reflection. The exact post-shock condition is applied to the left boundary and the bottom boundary from $x=-0.2$  to $x=0$. The reflecting boundary condition is applied to the wedge and upper boundary. Fig.\ref{fig:DMR_rho_contours} presents the density contours with the mesh size $h=1/80$ and $h=1/160$. The shock wave is captured sharply and the instability of the slip line is well resolved by the current scheme.
To further verify the subcell resolution of the current scheme, the density contours (blue) of the Mach stem near the wall are shown in Fig.\ref{fig:DMR_shock_thickness} where the black solid lines indicate main cells and the gray dashed lines indicate subcells. It can be clearly seen that the thickness of the Mach stem is smaller than the size of main cells for both $h=1/80$ and $h=1/160$. The shock wave spans across less than three subcells. Since discontinuity is allowed to exist on the interface between subcells, the resolution of shock waves can be much higher. 
Furthermore, the contours are smooth in the flow field. In contrast, for many other fourth-order methods, such as in Ref.\cite{GFu2017}, much finer mesh is usually required in order to achieve the same level of resolution, and it is very difficult to suppress the numerical noise occurring in the flow field, especially under the oblique shock wave.
These results demonstrates the capability of the current scheme to resolve hypersonic flows with strong robustness and high accuracy.

\subsection{Lid-driven cavity flow}
\begin{figure}[H]
  \centering
  \begin{varwidth}[t]{\textwidth}
  \vspace{0pt}
  \includegraphics[scale=0.3]{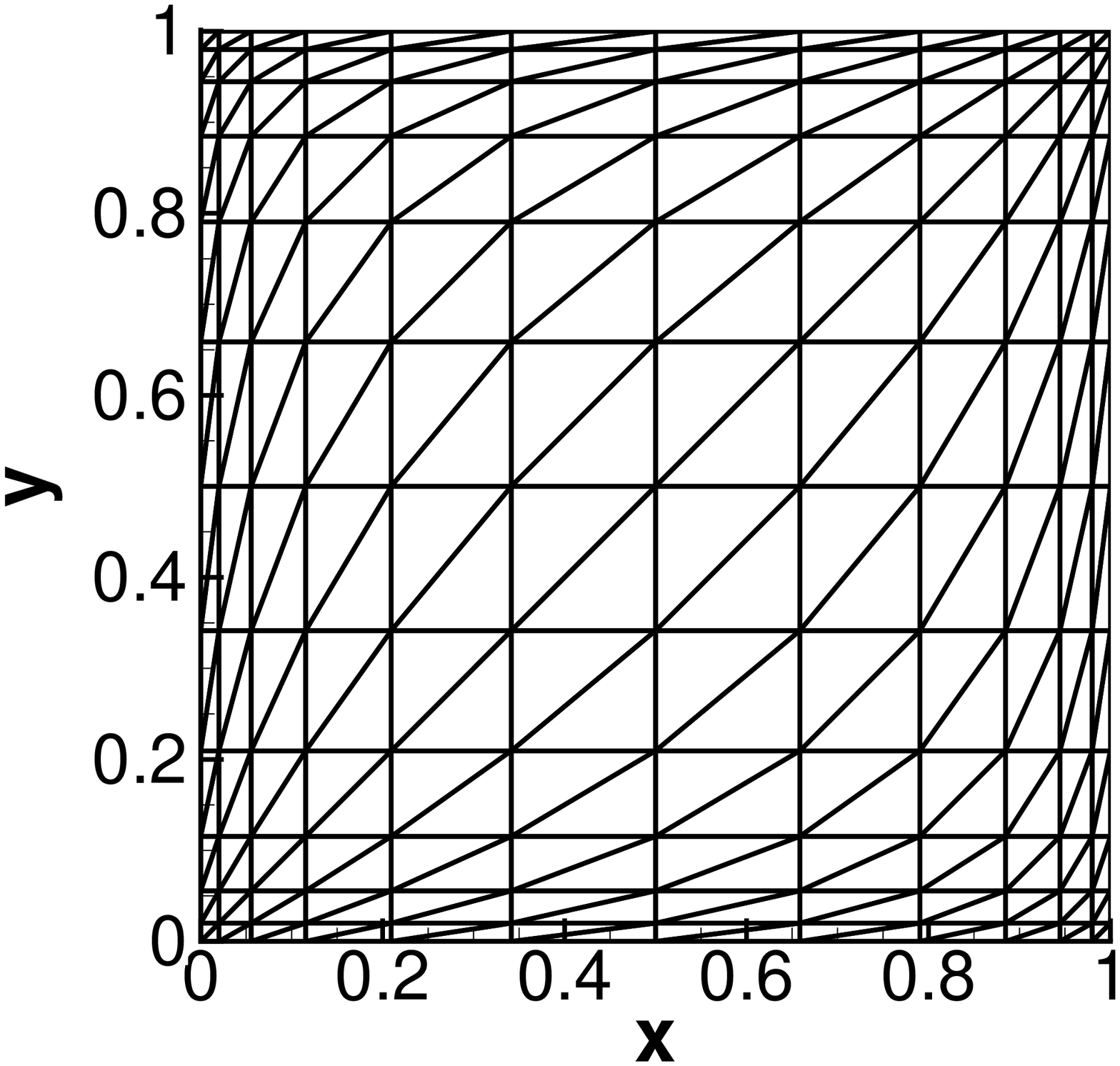}\\
  \end{varwidth}
  \qquad
  \begin{varwidth}[t]{\textwidth}
  \vspace{0pt}
  \includegraphics[scale=0.3]{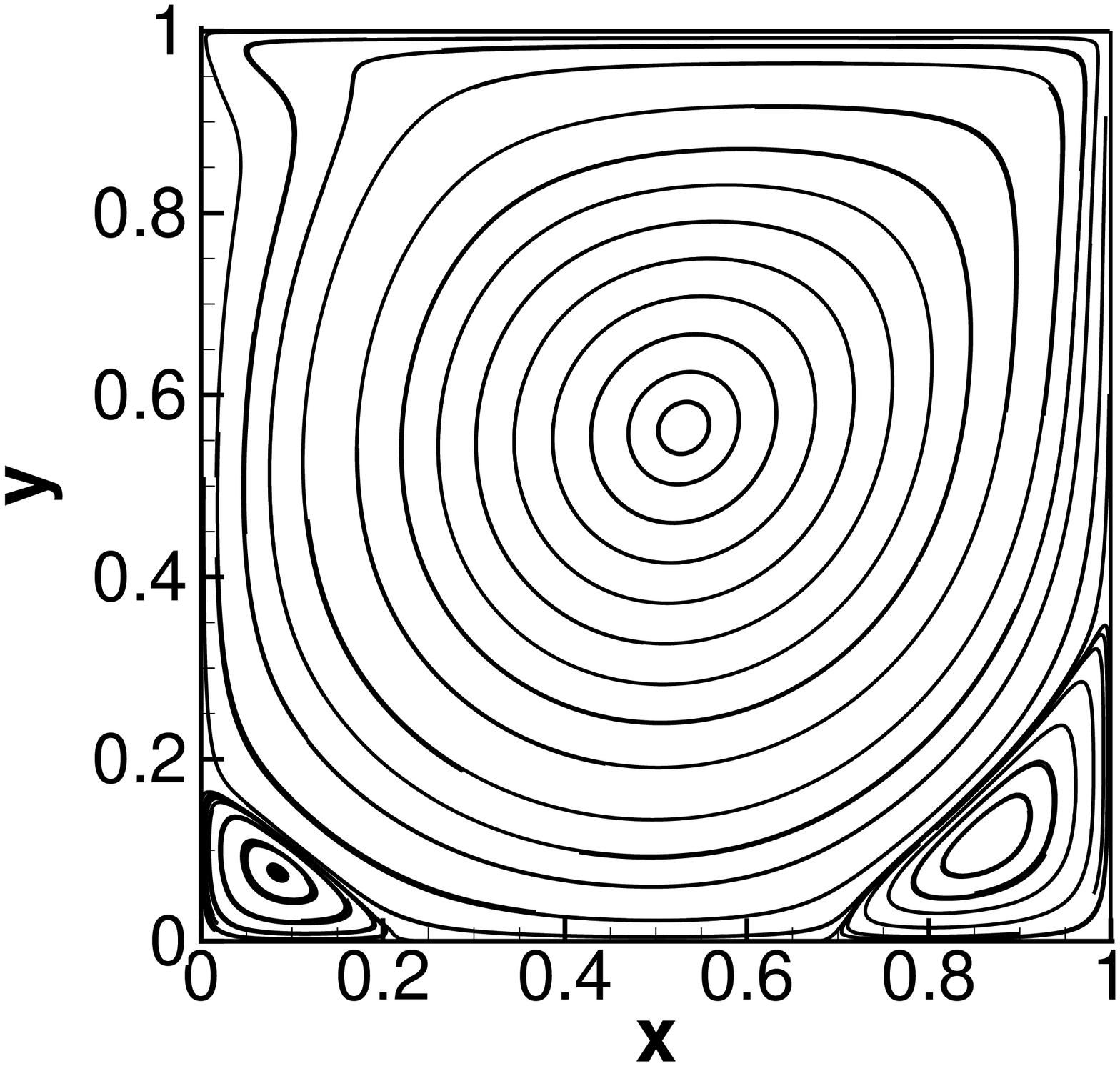}\\
  \end{varwidth}
  \captionsetup{labelsep=period}
  \caption{Computational mesh (left) and streamlines for Re=1000 (right) in lid-driven cavity flow.}
  \label{fig:Cavity_Mesh}
\end{figure}

The lid-driven cavity flow is a typical 2D incompressible viscous flow. It is bounded by a unit square $[0,1]\times[0,1]$. The upper wall is moving with the speed $U=1$, corresponding to a Mach number $\mathrm{Ma}=0.15$. Other walls are fixed. The non-slip and isothermal boundary conditions are applied to all boundaries with the temperature $T=1$. The initial flow is stationary with the density $\rho=1$ and temperature $T=1$. Two Reynolds numbers are considered, i.e., $\mathrm{Re}=1000$ and $\mathrm{Re}=3200$. As shown in Fig.\ref{fig:Cavity_Mesh}, the computational mesh contains $12\times12\times2$ elements with the minimum mesh size $h_{\min}=0.02$, and the stretching rate about $1.77$ near boundaries.

\begin{figure}[H]
  \centering
  \begin{varwidth}[t]{\textwidth}
  \vspace{0pt}
  \includegraphics[scale=0.33]{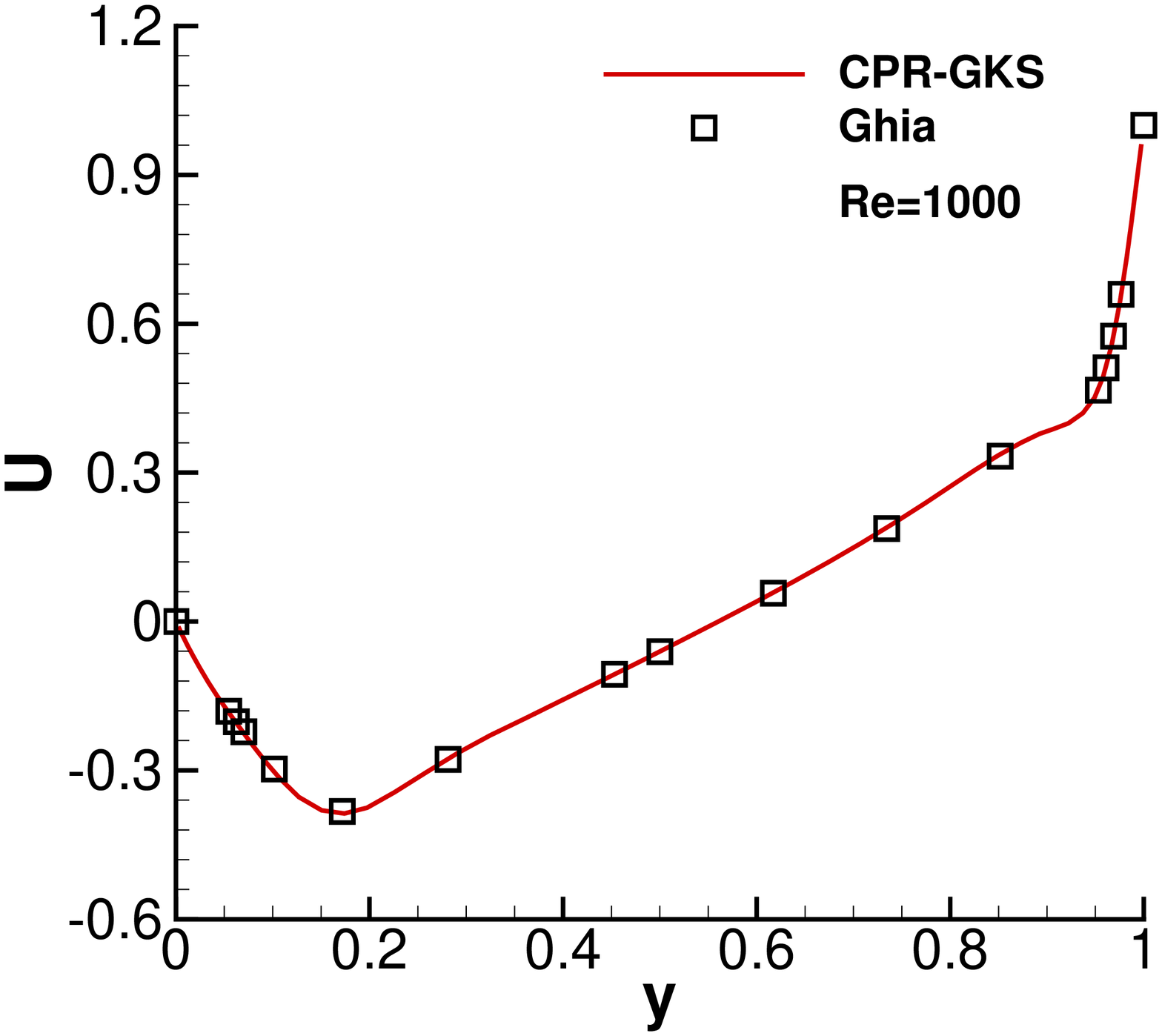}\\
  \end{varwidth}
  \qquad
  \begin{varwidth}[t]{\textwidth}
  \vspace{0pt}
  \includegraphics[scale=0.33]{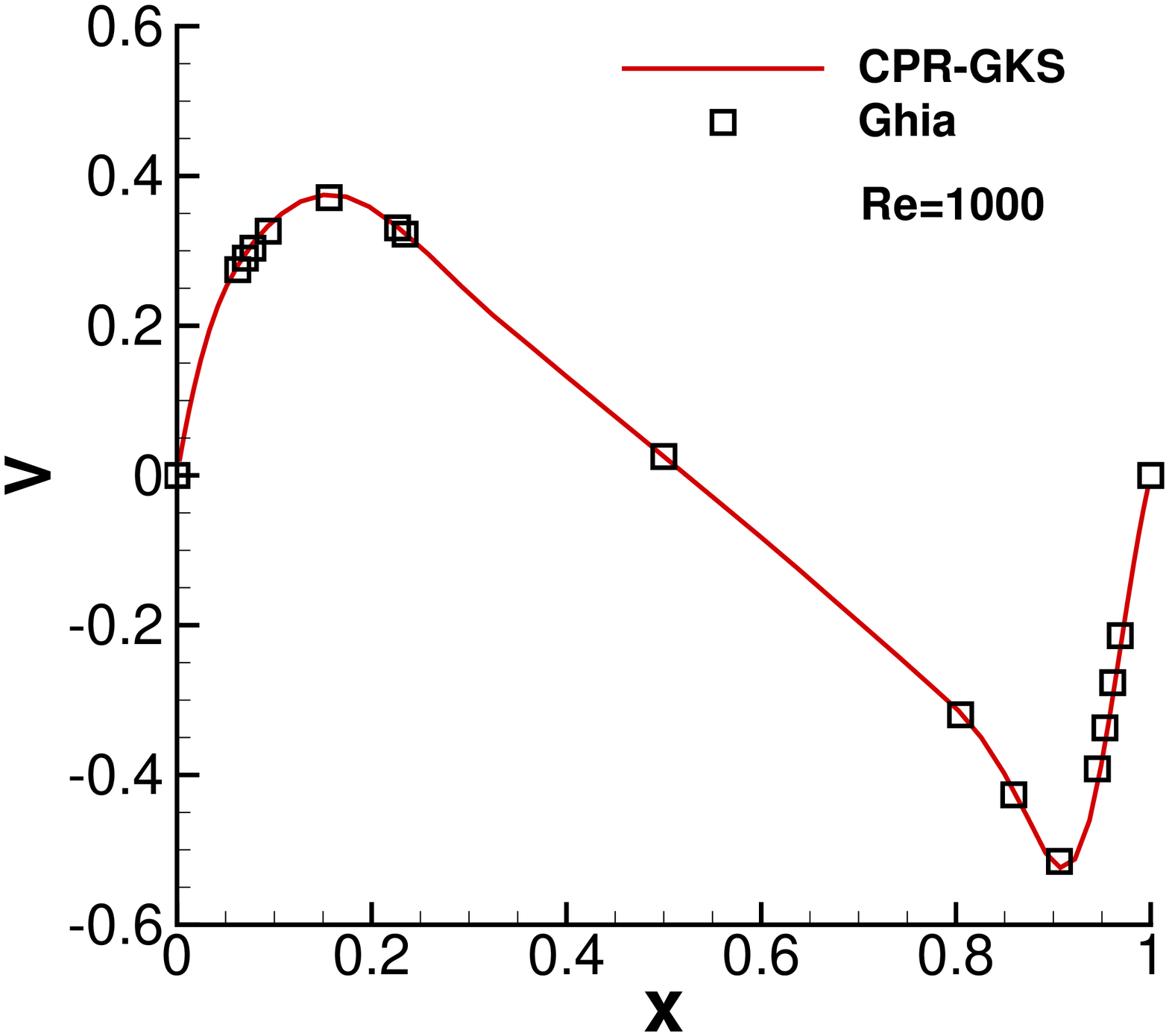}\\
  \end{varwidth}
  \qquad
  \captionsetup{labelsep=period}
  \caption{U-velocities along the vertical centerline and V-velocities along the horizontal centerline in lid-driven cavity flow with $\mathrm{Re}=1000$}
  \label{fig:Cavity_UV_Re1000}
\end{figure}

\begin{figure}[H]
  \centering
  \begin{varwidth}[t]{\textwidth}
  \vspace{0pt}
  \includegraphics[scale=0.33]{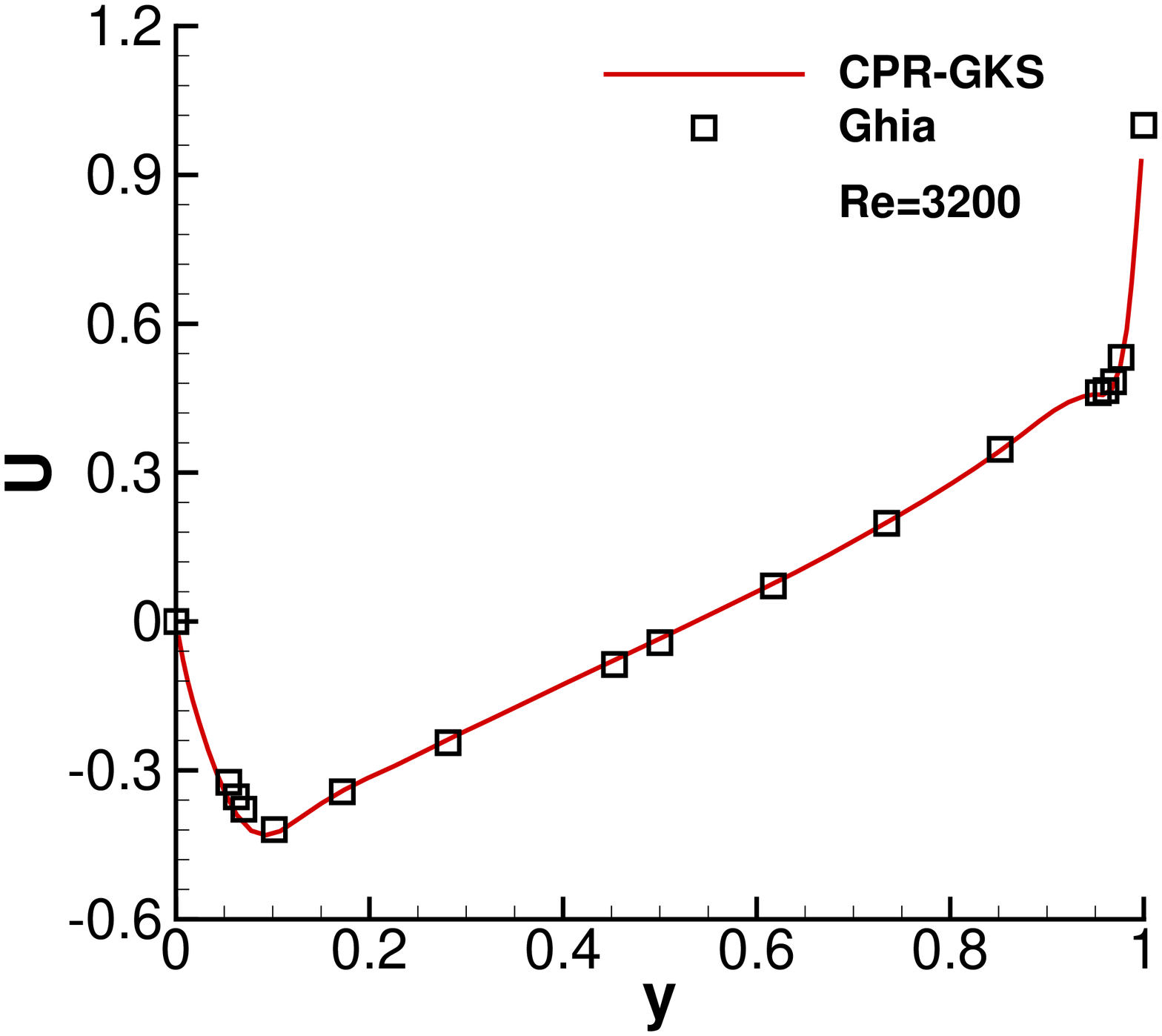}\\
  \end{varwidth}
  \qquad
  \begin{varwidth}[t]{\textwidth}
  \vspace{0pt}
  \includegraphics[scale=0.33]{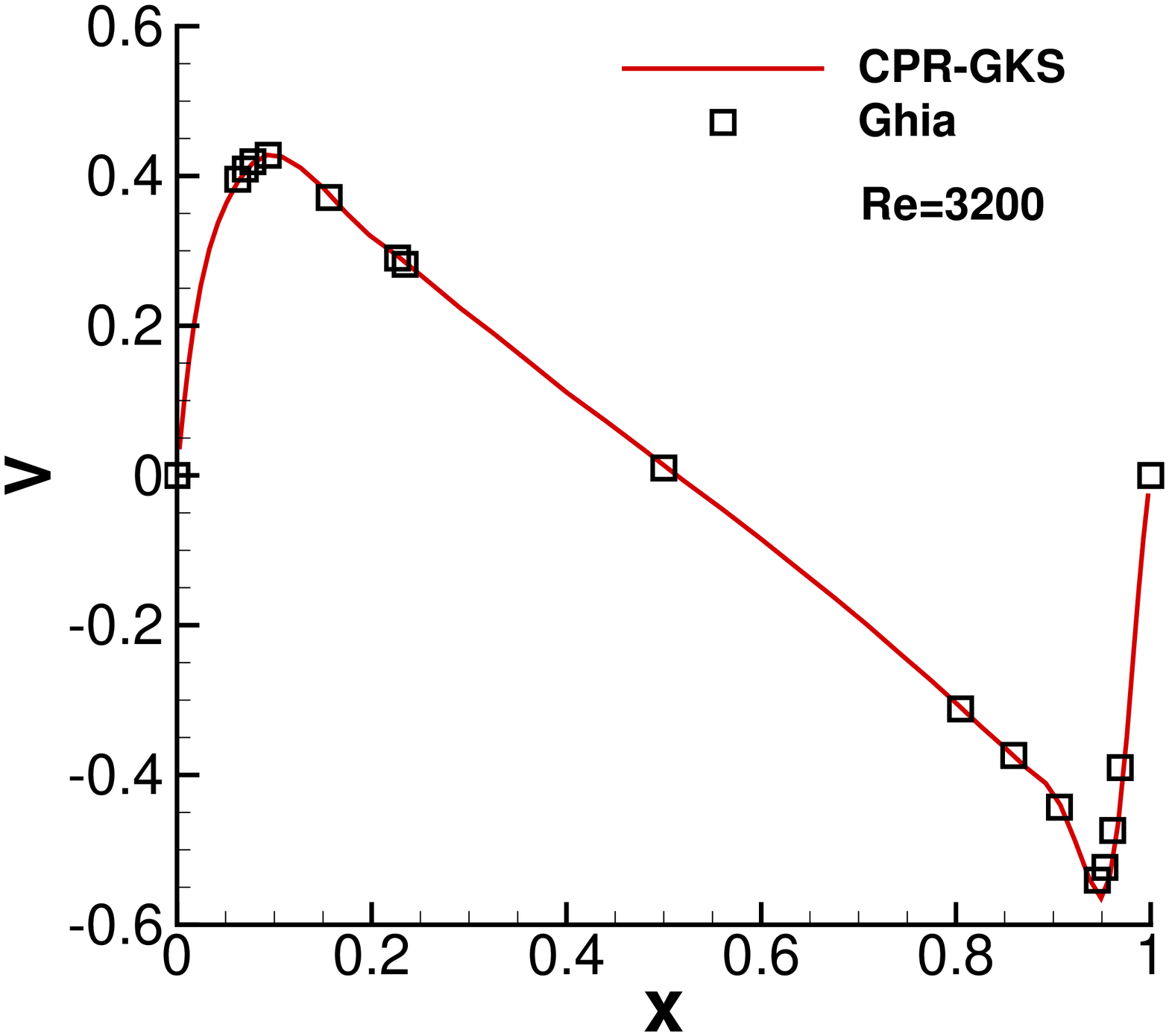}\\
  \end{varwidth}
  \captionsetup{labelsep=period}
  \caption{U-velocities along the vertical centerline and V-velocities along the horizontal centerline in lid-driven cavity flow with $\mathrm{Re}=3200$}
  \label{fig:Cavity_UV_Re3200}
\end{figure}

The streamlines for $\mathrm{Re}=1000$ are also shown in Fig.\ref{fig:Cavity_Mesh}. The primary flow structures including the primary and secondary vortices are well captured. Fig.\ref{fig:Cavity_UV_Re1000} and Fig.\ref{fig:Cavity_UV_Re3200} present the U-velocities along the vertical centerline and V-velocities along the horizontal centerline. It can be observed that, the velocity profiles agree very well with existing benchmark data for both $\mathrm{Re}=1000$ and $\mathrm{Re}=3200$. On such a coarse mesh with large stretching rate, it is challenging to resolve the velocity profiles accurately, especially the V-velocity profile for $\mathrm{Re}=3200$. The results demonstrate the high accuracy of the current scheme.

\subsection{Viscous shock tube problem}
To validate the performance of the current scheme in high-speed viscous flows, the viscous shock tube problem is simulated, which has been studied extensively \cite{Daru2004}. The flow is driven by a strong initial discontinuity at the center and bounded by a unit square. There exist complex unsteady interactions between the shock wave and the boundary layer, which requires not only strong robustness but also high resolution of a numerical scheme. The Reynolds number $\mathrm{Re}=200$ is chosen, which is based on  a constant dynamic viscosity $\mu=0.005$. The Prandtl number is $\mathrm{Pr}=0.73$. Considering the symmetry, the computational domain is set as $[0,1]\times[0,0.5]$ and the symmetrical condition is applied on the upper boundary. The non-slip and adiabatic conditions are adopted on other boundaries. The initial condition is
\begin{equation}\label{eq_VST_initial}
(\rho,U,V,p)=\begin{cases}
(120,0,0,120/\gamma),& 0\leq x \leq 0.5, \\
(1.2,0,0,1.2/\gamma),& 0.5\leq x \leq 1.
\end{cases}
\end{equation}

\begin{figure}[H]
  \centering
  \begin{varwidth}[t]{\textwidth}
  \vspace{0pt}
  \includegraphics[scale=0.3]{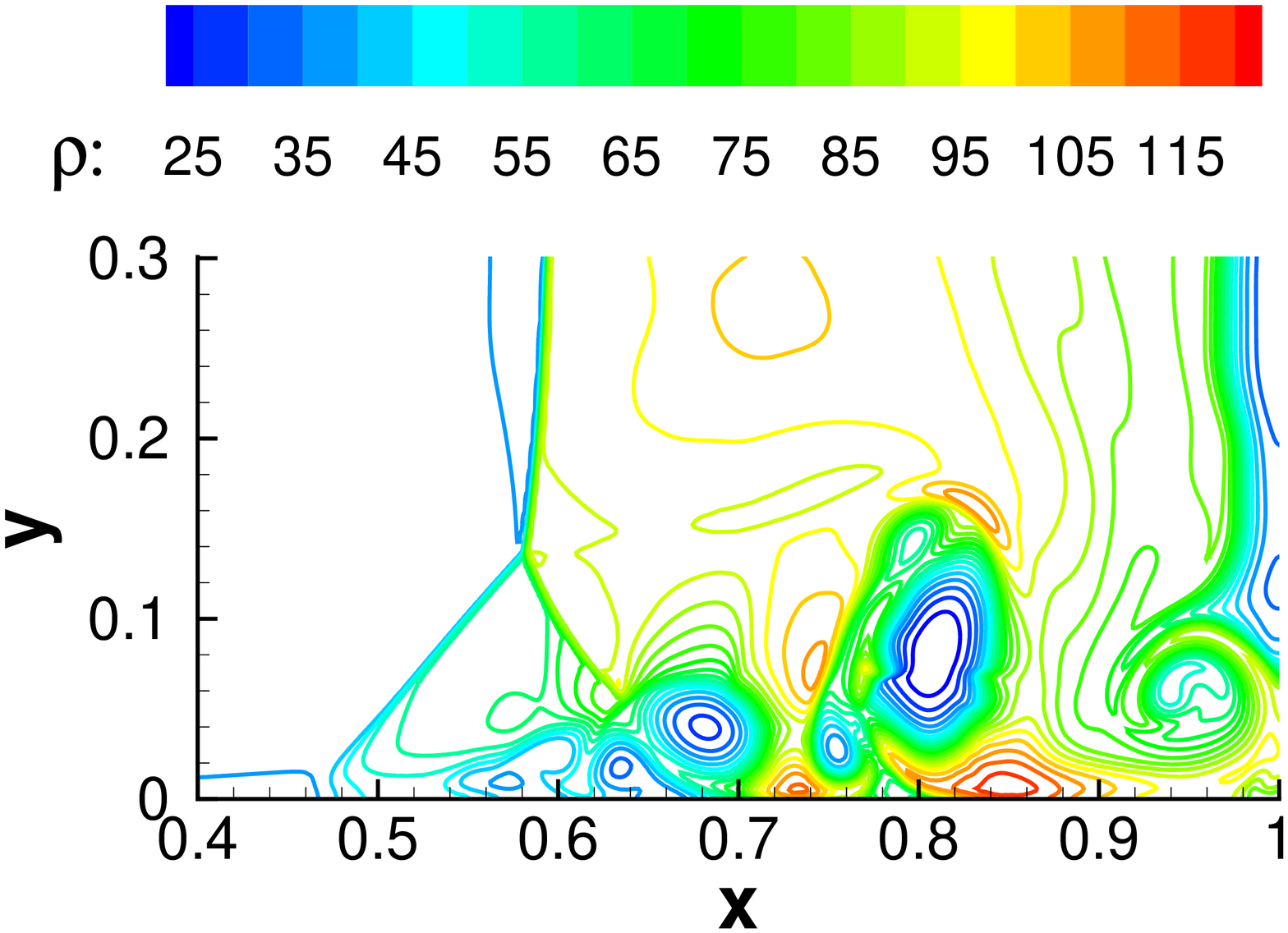}\\
  \end{varwidth}
  \qquad
  \begin{varwidth}[t]{\textwidth}
  \vspace{0pt}
  \includegraphics[scale=0.3]{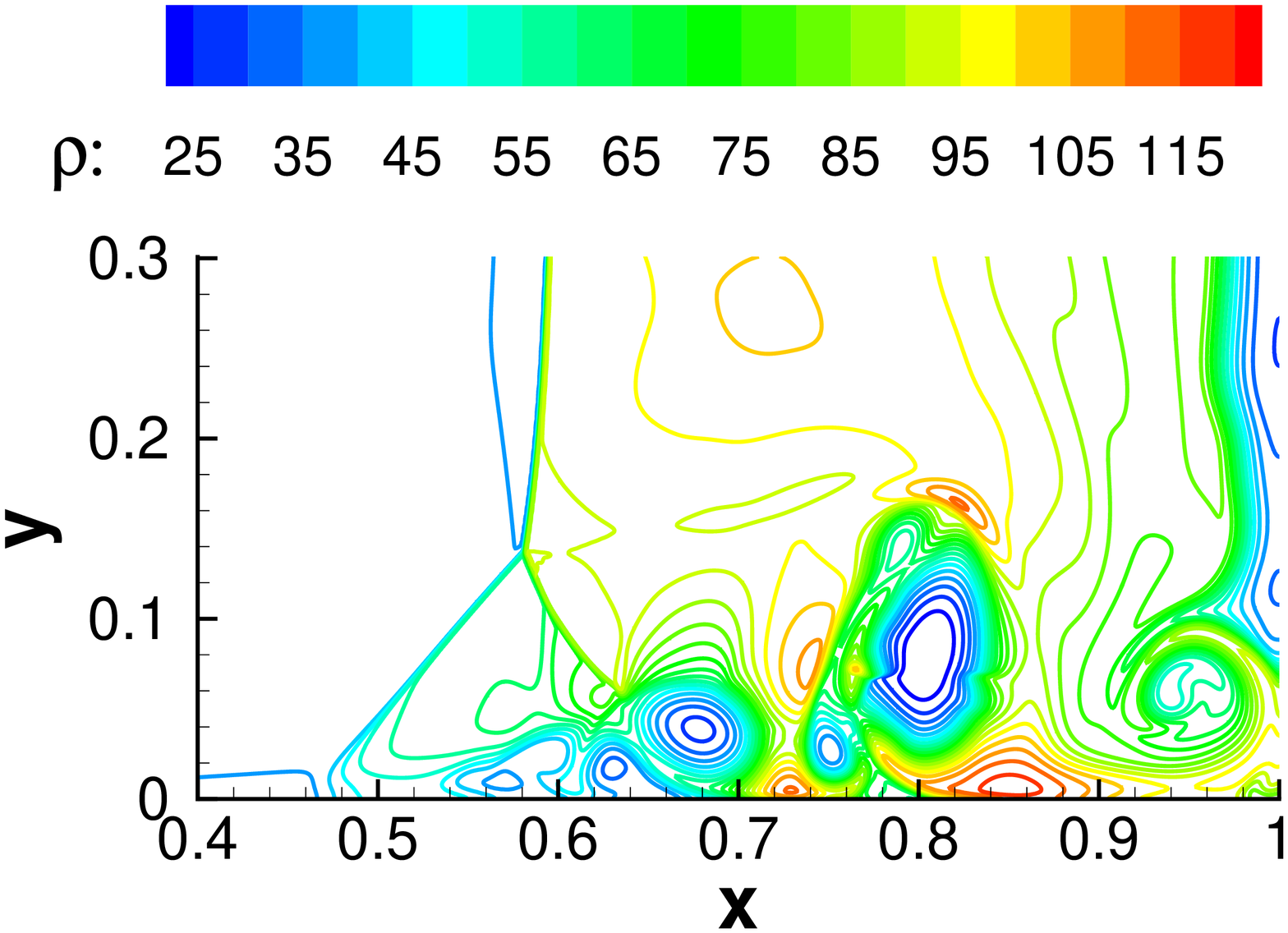}\\
  \end{varwidth}
  \qquad
  \captionsetup{labelsep=period}
  \caption{Density contours at t=1 in viscous shock tube flow with the mesh size $h=1/150$ (left) and $h=1/300$ (right). 20 uniform contours from 25 to 120.}
  \label{fig:VST_rho_2D}
\end{figure}

\begin{figure}[H]
  \centering
  \begin{varwidth}[t]{\textwidth}
  \vspace{0pt}
  \includegraphics[scale=0.35]{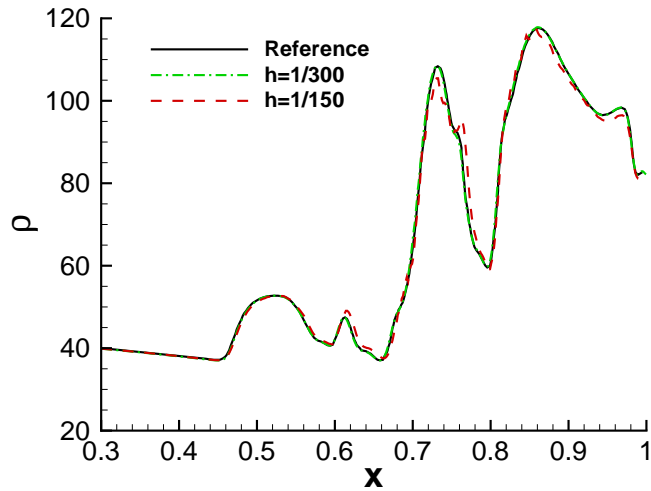}\\
  \end{varwidth}
  \captionsetup{labelsep=period}
  \caption{Density distribution along the bottom wall in the viscous shock tube flow.}
  \label{fig:VST_rho_wall}
\end{figure}

The result provided by a high-order GKS (HGKS) in Ref.\cite{GZZhou2018} is chosen as the reference data with the mesh size $h=1/1500$. Fig.\ref{fig:VST_rho_2D} shows the density contours at $t=1$ with the mesh size $h=1/150$ and $h=1/300$. It can be observed that the complex flow structure are well resolved by the current scheme, including the lambda shock and the vortex structures. Table~\ref{VST_Height} presents the height of primary vortex, achieving a good agreement with the reference data.
The density distribution along the bottom wall is presented in Fig.\ref{fig:VST_rho_wall}, which matches very well with the reference data, especially for the result with the mesh size $h=1/300$. Even with such a coarse mesh, the flow field can be well resolved by the current scheme, demonstrating the good performance of the current scheme in high-speed viscous flows.

\begin{table}[H]
\centering
\captionsetup{labelsep=period}
\caption{Comparison of the height of the primary vortex in viscous shock tube flow}
\label{VST_Height}
\renewcommand\arraystretch{1.3}
\begin{tabular}{lllll}
\hline
Scheme                       &  & $h$      &  & Height \\ \hline
\multirow{2}{*}{CPR-GKS}    &  & 1/150  &  & 0.168  \\
                             &  & 1/300  &  & 0.166  \\
HGKS                         &  & 1/1500 &  & 0.166  \\ \hline
\end{tabular}
\end{table}

\section{Conclusions}
A two-stage fourth-order gas-kinetic CPR method is developed for the compressible N-S equations on triangular meshes. It combines the efficient CPR framework with the robust and the efficient two-stage fourth-order temporal discretization, based on an efficient time-evolving gas-kinetic flux. Besides, a hybrid CPR/SCFV method is developed by extending a robust and accurate subcell finite volume limiting procedure  to the CPR framework to improve the resolution of discontinuities. Under the CPR framework, a fourth-order compact reconstruction is straightforward, simple and efficient, which avoids the difficulty of compactness encountered by traditional high-order finite volume GKS. Different from traditional CPR methods, the current method fully integrates the unique features of gas-kinetic flux solver. The inviscid and viscous fluxes are coupled and computed simultaneously. It is genuinely multi-dimensional by involving both normal and tangential variations in the gas distribution function. More importantly, with the time-evolving flux function, both flux and its first-order time derivative is available. The two-stage temporal discretization can therefore be extended to the CPR framework in a straightforward manner, which is more efficient than the multi-stage R-K method, by saving the computational cost of additional stages. In addition, with the help of SCFV limiting the subcell resolution of flow discontinuities is achieved, thus shock waves can be sharply and robustly captured even for extremely strong shock waves.

Several inviscid and viscous benchmark flows are simulated, ranging from nearly incompressible flows to supersonic flows with strong shock waves. It is verified that the current scheme achieves fourth-order accuracy in both space and time for the N-S equations, and the efficiency is higher than traditional CPR methods. For low-speed viscous flow, the current scheme has good agreement with benchmark data. For high-speed flows the smooth flow structures are accurately resolved while the strong shock waves can be well captured with subcell resolution. The present study demonstrate the high accuracy, efficiency and robustness of the current scheme, which makes it very competitive and promising in practical applications. 

Nevertheless, it should be noted that, the performance of the current hybrid CPR/SCFV method partly depends on the effectiveness of the shock detector. As only the second-order TVD limiter is considered in the SCFV method, if cells in smooth regions are mistakenly marked as troubled cells, the accuracy may be impaired. Thus, more investigations are required. To reduce the dependence on the shock detector, high-order limiter may also be implemented for the SCFV method. Besides,  
it is still necessary to make more comparisons with other typical high-order methods. The extension to three dimensions will be carried out in the future.

\section{Acknowledgments}
This work is supported by the National Natural Science Foundation of China (11672158, 91852109).
We also would like to acknowledge the technical support of PARATERA and the ``Explorer 100'' cluster system of Tsinghua National Laboratory for Information Science and Technology.

\section{References}

\bibliography{References}

\end{document}